\documentclass{elsart}
\usepackage{graphicx,psfrag}
\usepackage{amsmath}
\usepackage{amssymb}
\usepackage{calligra}
\begin{document}
\begin{frontmatter}

\title{An explicit solution for implicit time stepping in finite strain viscoelasticity.}
\author{A.V. Shutov\corauthref{cor}},
\corauth[cor]{Corresponding author. Tel.: +49-0-371-531-35024;
fax: +49-0-371-531-23419.}
\ead{alexey.shutov@mb.tu-chemnitz.de}
\author{R. Landgraf},
\author{J. Ihlemann}

\address{Chemnitz University of Technology,
Department of Solid Mechanics,
Str. d. Nationen 62, D-09111 Chemnitz, Germany}

\begin{abstract}

We consider the numerical treatment of one of the most popular finite strain models of the
viscoelastic Maxwell body. This model
is based on the multiplicative decomposition of
the deformation gradient,
combined with Neo-Hookean hyperelastic
relations between stresses and elastic strains.
The evolution equation is six dimensional.
For the corresponding local initial value problem, a
fully implicit integration procedure is considered,
and a simple explicit update formula is
derived. Thus, no local iterative procedure
is required, which makes the numerical
scheme more robust and efficient.
The resulting integration algorithm is
unconditionally stable and first order accurate.
The incompressibility constraint of
the inelastic flow is exactly preserved.
A rigorous proof of the
symmetry of the consistent tangent operator is provided. Moreover, some
properties of the numerical solution, like invariance under the change of the reference
configuration and positive energy dissipation within a time step, are discussed.
Numerical tests show that, in terms of accuracy, the proposed integration algorithm
is equivalent to the classical implicit scheme based
on the exponential mapping.
Finally, in order to check the stability of the algorithm numerically, a
representative initial boundary value problem
involving finite viscoelastic deformations is considered.
An FEM solution of the representative problem using MSC.MARC is presented.

\end{abstract}
\end{frontmatter}


\section{Introduction}

\begin{figure}[h]\centering
\psfrag{A}[m][][1][0]{$\varepsilon_{\text{i}}$}
\psfrag{B}[m][][1][0]{$\varepsilon_{\text{e}}$}
\psfrag{C}[m][][1][0]{$\varepsilon$}
\psfrag{D}[m][][1][0]{$\varepsilon_{\text{ii}}$}
\psfrag{E}[m][][1][0]{$\varepsilon_{\text{ie}}$}
\psfrag{F}[m][][1][0]{$a)$}
\psfrag{G}[m][][1][0]{$b)$}
\psfrag{H}[m][][1][0]{$c)$}
\scalebox{1.0}{\includegraphics{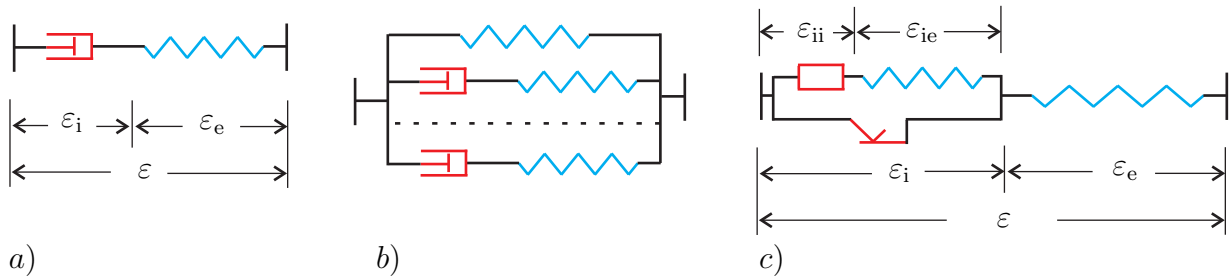}}
\caption{a) A one-dimensional Maxwell body consists of an elastic spring (Hooke body)
coupled in series with a viscous dashpot (Newton body),
b) Generalized Maxwell body, also known as Wiechert model (or Zener model in a special case) used
for better description of viscoelastic properties,
c) modified Schwedoff model used to represent nonlinear kinematic hardening.  \label{fig1}}
\end{figure}

Among idealized models of linear viscoelasticity, the
so-called Maxwell fluid (MF) is commonly encountered in material modelling \cite{Reiner, Haupt, Broec}.
The one-dimensional rheological interpretation of this model is shown in Fig. \ref{fig1}a.\footnote{A
two-dimensional rheological model of
the Maxwell fluid and its modifications can be
found in \cite{ShutovPaKr, ShutovIhle}.}
A series of Maxwell elements connected in parallel \cite{Wiechert} can be utilized
to represent viscoelastic properties of polymers
(Fig. \ref{fig1}b). In that case, the stresses acting in the Maxwell
elements can be associated with overstresses \cite{Haupt}.
Next, a slightly modified Maxwell element can be adopted to capture the
nonlinear kinematic hardening in metals (Fig. \ref{fig1}c). In that case, the corresponding Maxwell
stresses are interpreted as backstresses \cite{LionIJP, DettRes, Vladimirov, Feigen, ShutovKuprin}.
Moreover, within some phenomenological
approaches to metal plasticity, the distortional hardening
in metals can be captured using the modified MF \cite{ShutovPaKr, ShutovIhle}.
Other groups of materials like shape memory alloys \cite{Helm1, Helm3}
and biological tissues \cite{GassFor} can be modelled using MF.
Some applications of MF to finite deformations of geological structures \cite{SimoMeschke, PericCr} and
to fluid mechanics \cite{BalanTsakmakis} are reported in the literature as well.

In the finite strain range, numerous constitutive models of the MF exist
(see, among others, \cite{Leonov, JohnsonSegalman, LaMantia, PearsonMiddle,
NishiguchiA, NishiguchiB, Holzapfel, Drozdov, ReeseG, BalanTsakmakis, Haupt, HauptLi}).
Different variants were compared through numerical tests in \cite{GovRes, Landgraf}.
In this paper we consider one of the most popular models of the MF.
The corresponding constitutive equations are summarized in Section 2 of this work.
The model under consideration is a special case of the finite strain viscoplasticity model proposed by Simo and Miehe \cite{SimMieh},
and it has the same structure as the well known model of associative elastoplasticity considered by Simo \cite{Simo}.
These models were developed within the framework of multiplicative inelasticity in combination with
hyperelastic constitutive relations.
The corresponding inelastic flow rule\footnote{In this paper,
the evolution equation is referred to as ``inelastic flow rule''
in order to stress that the model is a special case of a viscoplastic model.} is six dimensional since
the inelastic spin plays no role due to elastic isotropy.
A version of the MF which is equivalent to the version of Simo \& Miehe was
considered later in material (Lagrangian) description by Lion \cite{LionAM}. This
Lagrangian formulation was adopted in \cite{LionHab, Helm1,
Hartmann2002, ShutovKrVisc, ShutovKrStab, ShutovPaKr}.
The spatial (Eulerian) constitutive equations proposed by
Simo \& Miehe were utilized later in the comprehensive study by Reese and Govindjee \cite{ReeseG},
and by many others (see, for instance, \cite{HubTsak, Nedjar, PericCr,
Kleuter, Hasanpour, HolLoug, Rauchs, Lejeunes}).

In the modern literature on numerics, much attention is paid to general
procedures which can be implemented to different types of constitutive
relations in a straightforward manner.
Because of their generality, such procedures are not always efficient
being compared to algorithms which make use of the special
structure of the underlying constitutive equations.
Due to the high prevalence of the Simo \& Miehe version of the Maxwell fluid in material modelling,
efficient and robust numerical integration of the
underlying evolution equations is a challenging task. \emph{The main purpose of this paper
is to report a new, simple, and efficient numerical procedure for this model}.

Since the corresponding initial value problem is typically stiff, implicit time
stepping methods should be implemented.
For the Simo \& Miehe version of the MF,
a discretized problem can be obtained using the operator split technique in combination
with exponential mapping and formulation in principal axes as
described in \cite{Simo} and \cite{ReeseG} (Eulerian approach). Alternatively,
the evolution equation formulated on the reference configuration
can be discretized as described in \cite{Hartmann} (Lagrangian approach).
In both cases, a system
of nonlinear algebraic equations is obtained, and
a local iterative procedure
is usually implemented to resolve the resulting
nonlinear problem (see, among others, \cite{ReeseG, Nedjar2,
Hartmann2002, HartmannHabil, PericCr, Helm2, Kleuter, ShutovKrKoo, Vladimirov, Hasanpour,
HolLoug, Rauchs, Lejeunes}).
Obviously, such iterative procedures can slow down the entire FEM simulation.
This problem may become especially important
if globally explicit FEM is considered.\footnote{In the case of explicit FEM,
the evaluation of the material routine at each point of Gauss integration
constitutes the major part of the overall computational effort.}
This publication is dealing with first-order accurate methods only.
For the discussion concerning the application of
higher order methods, the reader is referred
to \cite{ReeseHab, Hartmann2010, Eidel, EidelStumpf}.

In order to speed up the FEM computations,
much attention was paid to the construction of closed form solutions for
implicit schemes. For instance,
a simplified flow rule under the assumption of small elastic
strains was considered by Simo and Miehe \cite{SimMieh}
in order to get an explicit update formula for the
local implicit time stepping procedure.
For the same reason, another simplification
of the flow rule in case of small elastic strains was considered by
Reese and Govondjee \cite{ReeseG}. This simplified version
was implemented later in \cite{Johlitz}. Unfortunately, the simplifying assumption of
small elastic strains is not valid for many materials like
plastics, rubber, biological tissues etc. Moreover, if
the modified Maxwell body is used to capture nonlinear kinematic hardening in metals,
a general finite strain version of the model must be utilized
as well.\footnote{In fact, although the elastic strains in metals are typically
small ($\varepsilon_{\text{e}} \rightarrow 0$ in Fig. \ref{fig1}c), the
conservative part ($\varepsilon_\text{ie}$ in Fig. \ref{fig1}c) of the
inelastic strain may become finite.}
Another approach to closed form solution is based on special assumptions concerning the energy storage.
In particular, a quadratic logarithmic strain energy (so-called Hencky strain energy) can be assumed in order
to simplify the numerical treatment of the material model \cite{Meng}.
Unfortunately, this assumption would yield unrealistic results in case
of large elastic strains. Thus, again,
the applicability area is limited to moderate elastic strains.
In this work, a simple explicit update formula is presented for
the original finite strain version with Neo-Hookean potential.
Interestingly, this explicit solution for the general case is even
more compact and simple than the solutions presented in \cite{SimMieh} and \cite{ReeseG}
for the special case of small elastic strains or the
solution in \cite{Meng} for quadratic logarithmic strain energy.
For the new method, the computational effort per single time step is even smaller than the
effort required within the explicit time stepping.

The inelastic flow is assumed to
be incompressible, and the algorithm presented in this work
preserves this incompressibility constraint.
A classical model of finite strain viscoplasticity which contains
the Simo \& Miehe version of the MF was considered in  \cite{ShutovKrStab}.
As it was shown in \cite{ShutovKrStab},
the exact solution to the initial value
problem is exponentially stable with respect to small perturbations
of the initial data, if the incompressibility constraint is not violated.
For such material models, the numerical schemes which exactly
preserve the incompressibility are advantageous
due to the \emph{suppressed error accumulation} \cite{ShutovKrStab}.
This theoretical result is confirmed by numerical tests
presented in the current paper.

Dealing with the constitutive equations written in Lagrangian form,
it can be shown that they are invariant
under the isochoric change of the reference configuration \cite{ShutovPfeIhl}.
The same invariance property
can be formulated for the numerical solution as well.
Obviously, the numerical algorithms which exactly retain this invariance property
are advantageous. In this work, it is proved that the advocated algorithm
retains the invariance of the solution.

We conclude this introduction with a few words regarding notation.
Throughout this article, bold-faced symbols denote first- and second-rank tensors in $\mathbb{R}^3$.
A coordinate-free tensor formalism is used in this work \cite{Itskov, ShutovKrKoo}.
In this work, $\mathbf{1}$ stands for the second-rank identity tensor.
The deviatoric part of a tensor is defined as
$\mathbf A^{\text{D}} := \mathbf A - \frac{1}{3} \text{tr}(\mathbf A) \mathbf 1$, where
$\text{tr}(\mathbf A)$ stands for the trace.
The material time derivative is denoted by dot: $\frac{d}{d t} \mathbf A = \dot{\mathbf A}$.
The overline $\overline{(\cdot)}$ denotes the unimodular
part of a tensor such that $\overline{\mathbf{A}}=(\det \mathbf{A})^{-1/3} \mathbf{A}$.
The inverse of transposed tensor is denoted by $ \mathbf{A}^{-\text{T}}$.
The positive definiteness of a tensor $\mathbf{A}$ is symbolically denoted by $\mathbf{A} > 0$.

\section{System of constitutive equations}

\subsection{Lagrangian formulation}

Let us consider a finite strain model of Maxwell fluid. This model is covered as a
special case by the viscoplasticity model presented by Simo and Miehe \cite{SimMieh}.
The Lagrangian formulation of the model follows the presentation of Lion \cite{LionAM}.
We start with the multiplicative decomposition of the deformation gradient $\mathbf F$
into the elastic part $\hat{\mathbf F}_{\text{e}}$ and the inelastic part $\mathbf F_{\text{i}}$
\begin{equation*}\label{split1}
\mathbf F = \hat{\mathbf F}_{\text{e}} \mathbf F_{\text{i}}.
\end{equation*}
Along with the well-known
right Cauchy-Green tensor
${\mathbf C}=\mathbf F^{\text{T}} \mathbf F$, we introduce
the inelastic right Cauchy-Green tensor as
\begin{equation*}\label{intIntv}
{\mathbf C}_{\text{i}}=\mathbf F^{\text{T}}_{\text{i}}
\mathbf F_{\text{i}}.
\end{equation*}
The elastic right Cauchy-Green tensor $\hat{\mathbf C}_{\text{e}}$
and the elastic Green tensor $\hat{\mathbf \Gamma}_{\text{e}}$ are defined by
\begin{equation*}\label{cecie}
\hat{\mathbf C}_{\text{e}}  :=
\hat{\mathbf F}_{\text{e}}^{\text{T}} \hat{\mathbf F}_{\text{e}}, \quad
\hat{\mathbf \Gamma}_{\text{e}} = \frac{1}{2} (\hat{\mathbf C}_{\text{e}} - \mathbf{1}).
\end{equation*}
Next, we introduce the inelastic velocity gradient
$\hat{\mathbf L}_{\text{i}} $ and the covariant Oldroyd derivative (with respect to the
intermediate configuration)
\begin{equation*}\label{altol}
\hat{\mathbf L}_{\text{i}} = \dot{\mathbf F}_{\text{i}} \mathbf F^{-1}_{\text{i}}, \quad
\stackrel{\triangle} {(\cdot)} = \frac{d}{d t}(\cdot) +
\hat{\mathbf L}^{\text{T}}_{\text{i}}
(\cdot) + (\cdot) \hat{\mathbf L}_{\text{i}}.
\end{equation*}
The inelastic Almansi strain tensor $\hat{\mathbf \Gamma}_{\text{i}} $ and the inelastic
strain rate tensor $\hat{\mathbf D}_{\text{i}} $ are defined through
 \begin{equation*}\label{defgae}
\hat{\mathbf \Gamma}_{\text{i}} := \frac{\displaystyle 1}{\displaystyle 2}
(\mathbf 1 - \mathbf F_{\text{i}}^{-\text{T}} \mathbf F_{\text{i}}^{-1}), \quad
\hat{\mathbf D}_{\text{i}} = \frac{1}{2}
(\hat{\mathbf L}_{\text{i}} + \hat{\mathbf L}^{\text{T}}_{\text{i}} ).
\end{equation*}
After some straightforward computations (cf. \cite{Haupt}), one gets
\begin{equation}\label{DOldrDer}
\hat{\mathbf D}_{\text{i}} = \stackrel{\triangle}{\hat{\mathbf \Gamma}}_{\text{i}}.
\end{equation}
Let $\mathbf T$ be the Cauchy stress tensor.
The Kirchhoff stress tensor $\mathbf S$,
the 2nd Piola-Kirchhoff tensor $\hat{\mathbf S}$ operating on the intermediate configuration
and the 2nd Piola-Kirchhoff tensor $\tilde{\mathbf T}$ operating on the reference configuration are defined through
\begin{equation}\label{Kirch}
\mathbf S  : = (\text{det} \mathbf F) \mathbf T, \quad
\hat{\mathbf S} : = {\hat{\mathbf F}_{\text{e}}^{-1}} \mathbf S {\hat{\mathbf F}_{\text{e}}^{-\text{T}}}, \quad
\tilde{\mathbf T} : = {\mathbf F^{-1}} \mathbf S {\mathbf F^{-\text{T}}}.
\end{equation}
Let $\psi$ be the free energy per unit mass. It is postulated in the Neo-Hookean form as
\begin{equation*}\label{spec1}
\rho_{\scriptscriptstyle \text{R}}  \psi(\hat{\mathbf{C}}_{\text{e}})=
\frac{\mu}{2} \big( \text{tr} \overline{\hat{\mathbf{C}}_{\text{e}}} - 3 \big).
\end{equation*}
Here, $\rho_{\scriptscriptstyle \text{R}} > 0$ stands for the mass density with respect to the reference
configuration; $\mu \geq 0$ is the shear modulus. Next, noting
that $\hat{\mathbf C}_{\text{e}} = \mathbf{1}+2 \hat{\mathbf{\Gamma}}_{\text{e}}$,
a hyperelastic stress-strain relation is introduced
on the intermediate configuration:
\begin{equation}\label{potent}
\hat{\mathbf S} = \rho_{\scriptscriptstyle \text{R}}
\frac{\displaystyle \partial \psi(\mathbf{1}+2 \hat{\mathbf{\Gamma}}_{\text{e}})}
{\displaystyle \partial \hat{\mathbf{\Gamma}}_{\text{e}}}.
\end{equation}

Isothermal processes are considered in this study. The Clausius-Duhem inequality requires that
the specific internal dissipation $\delta_{\text{i}}$ remains non-negative (see \cite{Haupt})
\begin{equation}\label{cld}
\delta_{\text{i}} := \frac{1}{\rho_{\scriptscriptstyle \text{R}}} \tilde{\mathbf T} : \dot{\mathbf E} - \dot{\psi} \geq 0,
\end{equation}
where $\mathbf E : = \frac{\displaystyle 1}{\displaystyle 2} (\mathbf C - \mathbf 1)$ stands for the
Green strain tensor.
Using \eqref{potent} and taking the isotropy of the free energy function into account,
this inequality is reduced to
\begin{equation}\label{cld2}
\rho_{\scriptscriptstyle \text{R}} \delta_{\text{i}} =
(\hat{\mathbf{C}}_{\text{e}}  \hat{\mathbf S}) : \stackrel{\triangle} {\hat{\mathbf{\Gamma}}}_{\text{i}}  \geq 0.
\end{equation}

An evolution equation is postulated so that
inequality \eqref{cld2} holds for arbitrary
mechanical loadings (cf. \cite{LionAM})
\begin{equation}\label{evol}
\stackrel{\triangle}{\hat{\mathbf{\Gamma}}}_{\text{i}} = \frac{1}{2 \eta}
(\hat{\mathbf{C}}_{\text{e}}  \hat{\mathbf S})^{\text{D}},
\end{equation}
where $\eta \geq 0$ is a material parameter (Newtonian viscosity).\footnote{In some
works, the viscosity $\eta$ is replaced by the parameter $\varkappa >0$ such that
$\varkappa = 1/\eta$.
Furthermore, the relaxation time can be introduced by $\tau = \eta / \mu$.
The case of process-dependent viscosity was considered,
among others, in \cite{Leonov, LaMantia, LionAM, ShutovKuprin, Koprowski}.}
In view of \eqref{DOldrDer}, an equivalent formulation of this flow
rule is given by
\begin{equation}\label{evol2}
\hat{\mathbf{D}}_{\text{i}} = \frac{1}{2 \eta}
(\hat{\mathbf{C}}_{\text{e}}  \hat{\mathbf S})^{\text{D}}.
\end{equation}
Note that this flow rule is six dimensional since
the rotational part
$\text{skew}(\hat{\mathbf L}_{\text{i}})$
 drops out from the constitutive
relations due to the elastic isotropy.
Moreover, since $\text{tr} (\stackrel{\triangle}{\hat{\mathbf{\Gamma}}}_{\text{i}}) =
\text{tr} (\hat{\mathbf{D}}_{\text{i}}) =0$,
the inelastic flow described by \eqref{evol} in incompressible.

In order to simplify the numerical treatment of the model, the
constitutive equations can be transformed to the reference configuration.
In particular, the free energy takes the form
\begin{equation}\label{freeen2}
\psi=\psi(\mathbf C {\mathbf C_{\text{i}}}^{-1})=
\frac{\mu}{2 \rho_{\scriptscriptstyle \text{R}} } \big(
\text{tr} \overline{\mathbf C {\mathbf C_{\text{i}}}^{-1}} - 3 \big).
\end{equation}
Using \eqref{Kirch} and \eqref{potent}, one gets for the 2nd Piola-Kirchhoff stress tensor
\begin{equation*}\label{2PKddd}
\tilde{\mathbf T}  =
2 \rho_{\scriptscriptstyle \text{R}}
\frac{\displaystyle \partial \psi(\mathbf C {\mathbf C_{\text{i}}}^{-1})}
{\displaystyle \partial \mathbf{C}}\big|_{\mathbf C_{\text{i}} =
\text{const}}.
\end{equation*}
Substituting \eqref{freeen2} into this relation, one gets
\begin{equation}\label{2PKd}
\tilde{\mathbf T}  = \mu \ \mathbf C^{-1} (\overline{\mathbf C} \mathbf C_{\text{i}}^{-1})^{\text{D}}.
\end{equation}
Next, we note that
\begin{equation*}\label{trace}
\text{tr} (\hat{\mathbf{C}}_{\text{e}}  \hat{\mathbf S}) =
\text{tr} (\mathbf C \tilde{\mathbf T}).
\end{equation*}
Thus, the pull-back of the deviatoric part of the Mandel tensor $\hat{\mathbf{C}}_{\text{e}}  \hat{\mathbf S}$ is given by
\begin{equation*}\label{puba2}
{\mathbf F}_{\text{i}}^{\text{T}} (\hat{\mathbf{C}}_{\text{e}}  \hat{\mathbf S})^{\text{D}} {\mathbf F}_{\text{i}} =
\mathbf C \tilde{\mathbf T} {\mathbf C}_{\text{i}}  - \frac{1}{3} \text{tr} (\hat{\mathbf{C}}_{\text{e}}  \hat{\mathbf S})
\ {\mathbf C}_{\text{i}} \stackrel{\eqref{trace}}{=}
\big( \mathbf C \tilde{\mathbf T} \big)^{\text{D}} \mathbf C_{\text{i}}.
\end{equation*}
Using this relation in combination with the evolution equation \eqref{evol}, we get
\begin{equation*}\label{puba4}
\dot{\mathbf C}_{\text{i}} = 2 {\mathbf F}_{\text{i}}^{\text{T}}
\stackrel{\triangle}{\hat{\mathbf{\Gamma}}}_{\text{i}} {\mathbf F}_{\text{i}} \stackrel{\eqref{evol}}{=}
 \frac{1}{\eta} {\mathbf F}_{\text{i}}^{\text{T}} (\hat{\mathbf{C}}_{\text{e}}  \hat{\mathbf S})^{\text{D}} {\mathbf F}_{\text{i}}
  \stackrel{\eqref{puba2}}{=}  \frac{1}{\eta} \big( \mathbf C \tilde{\mathbf T} \big)^{\text{D}} \mathbf C_{\text{i}}.
\end{equation*}
Taking \eqref{2PKd} into account, we have
\begin{equation}\label{puba44}
\dot{\mathbf C}_{\text{i}} = \frac{\mu}{\eta}  \big( \overline{\mathbf C}
\mathbf C^{-1}_{\text{i}}  \big)^{\text{D}} \mathbf C_{\text{i}}.
\end{equation}
The system of constitutive equations \eqref{2PKd} and \eqref{puba44} is closed by specifying initial condinions
\begin{equation}\label{initCond}
\mathbf C_{\text{i}}|_{t=t^0} = \mathbf C_{\text{i}}^0.
\end{equation}

Note that the exact solution to the evolution equation \eqref{puba44} has
the following geometric property
\begin{equation}\label{geopro}
\mathbf{C}_{\text{i}}(t) \in \mathbb{M} \quad \text{if} \quad
\mathbf{C}^0_{\text{i}} \in \mathbb{M},
\end{equation}
where the manifold $\mathbb{M}$ is a set of symmetric unimodular tensors
\begin{equation*}\label{UnimodulM}
\mathbb{M} := \big\{ \mathbf A \in Sym: \text{det} \mathbf A =1 \big\}.
\end{equation*}
It follows from \eqref{geopro} that $\mathbf{C}_{\text{i}}$ remains positive definite if $\mathbf{C}^0_{\text{i}} > 0$.\footnote{
The eigenvalues of $\mathbf{C}_{\text{i}}$ are continuous functions of time. Therefore,
if all of the eigenvalues were positive at some time instance
and their product remains constant in time, then they remain positive.}

\textbf{Remark 1} Interestingly,
an evolution equation in the form \eqref{puba44} was presented for an
alternative model of the MF by other authors (see equation (10.147) in \cite{Haupt}),
although the flow rule considered on the intermediate configuration
is given by $\hat{\mathbf{D}}_{\text{i}} = \frac{1}{\eta} \hat{\mathbf S} $, which
differs essentially from the flow rule \eqref{evol2}.
Since both evolution equations coincide,
the approach advocated in this paper can
be applied to the version presented in \cite{Haupt}
without any modifications whatsoever.

\subsection{Eulerian formulation}

Now let us check that the model presented in this section
coincides with the model of Simo and Miehe \cite{SimMieh}, formulated
within the Eulerian approach. Here, elastic isotropy is considered.
First,  note that
\begin{equation}\label{SimMie}
2 {\hat{\mathbf{D}}}_{\text{i}} = {\mathbf F}_{\text{i}}^{-\text{T}}
\dot{\mathbf C}_{\text{i}} {\mathbf F}_{\text{i}}^{-1}.
\end{equation}
Due to the elastic isotropy, according to the
evolution equation \eqref{evol2}, ${\hat{\mathbf{D}}}_{\text{i}}$ commutes with $\hat{\mathbf{C}}_{\text{e}}$.
In particular, one gets from \eqref{evol2}
\begin{equation}\label{SimMie2}
\hat{\mathbf{C}}_{\text{e}} 2 \hat{\mathbf{D}}_{\text{i}}
\hat{\mathbf{C}}_{\text{e}}^{-1} = \frac{1}{\eta}
(\hat{\mathbf{C}}_{\text{e}}  \hat{\mathbf S})^{\text{D}}.
\end{equation}
Substituting \eqref{SimMie} into \eqref{SimMie2} and taking into account
that $\frac{d}{d t} ({\mathbf C}_{\text{i}}^{-1} ) = -
{\mathbf C}_{\text{i}}^{-1} \dot{{\mathbf C}}_{\text{i}} {\mathbf C}_{\text{i}}^{-1}$ one
gets
\begin{equation}\label{SimMie3}
-\hat{\mathbf{F}}^{\text{T}}_{\text{e}} \mathbf{F}
\frac{d}{d t} ({\mathbf C}_{\text{i}}^{-1} ) \mathbf{F}^{\text{T}}
\hat{\mathbf{F}}_{\text{e}}^{-\text{T}}
\hat{\mathbf{F}}_{\text{e}}^{-1} \hat{\mathbf{F}}_{\text{e}}^{-\text{T}}   = \frac{1}{\eta}
(\hat{\mathbf{C}}_{\text{e}}  \hat{\mathbf S})^{\text{D}}.
\end{equation}
For what follows, we introduce the elastic left Cauchy-Green tensor
$\hat{\mathbf{B}}_{\text{e}} = \hat{\mathbf{F}}_{\text{e}} \hat{\mathbf{F}}_{\text{e}}^{\text{T}}$, and
note that
\begin{equation*}\label{SimMie4}
\mathbf{S}^{\text{D}} = \hat{\mathbf{F}}_{\text{e}}^{-\text{T}}
(\hat{\mathbf{C}}_{\text{e}}  \hat{\mathbf S})^{\text{D}}\hat{\mathbf{F}}_{\text{e}}^{\text{T}}.
\end{equation*}
Multiplying \eqref{SimMie3} with $\hat{\mathbf{F}}^{-\text{T}}_{\text{e}}$ from the left and
with $\hat{\mathbf{F}}^{\text{T}}_{\text{e}}$ from the right one gets using \eqref{SimMie4}
\begin{equation}\label{SimMie5}
-\mathbf{F}
\frac{d}{d t} ({\mathbf C}_{\text{i}}^{-1} ) \mathbf{F}^{\text{T}}
\hat{\mathbf{B}}_{\text{e}}^{-1}
 = \frac{1}{\eta} \mathbf S^{\text{D}}.
\end{equation}
By introducing the covariant Oldroyd rate\footnote{In order to sress that this Oldroyd
rate is also a Lee derivative, it can be denoted
by $\text{\calligra{L}}_{v} $.}
\begin{equation*}\label{OldrCovar}
\mathfrak{O} (\mathbf{A}) =
\text{\calligra{L}}_{v}  (\mathbf{A}) : =
\mathbf{F} \frac{d}{dt} (\mathbf{F}^{-1} \mathbf{A} \mathbf{F}^{-\text{T}})  \mathbf{F}^{\text{T}}
= \dot{\mathbf{A}} -  \mathbf{L}  \mathbf{A} -  \mathbf{A}  \mathbf{L}^{\text{T}},
\end{equation*}
one gets
\begin{equation*}\label{OldrCovar2}
\mathfrak{O} (\hat{\mathbf{B}}_{\text{e}}) =
\text{\calligra{L}}_{v}  (\hat{\mathbf{B}}_{\text{e}}) = \mathbf{F}
\frac{d}{d t} ({\mathbf C}_{\text{i}}^{-1} ) \mathbf{F}^{\text{T}}.
\end{equation*}
Thus, the flow rule \eqref{SimMie5} takes the well-known form
\begin{equation}\label{SimMie82}
-\text{\calligra{L}}_{v}  (\hat{\mathbf{B}}_{\text{e}})
\hat{\mathbf{B}}_{\text{e}}^{-1}
 = \frac{1}{\eta} \mathbf S^{\text{D}}.
\end{equation}
This equation was covered as a special case by Simo and Miehe \cite{SimMieh}
(see equations (2.19a) and (2.26) in \cite{SimMieh})
as well as by Reese and Govindjee \cite{ReeseG}.
In case of the Neo-Hookean potential \eqref{freeen2}, we have
$\mathbf S = \mathbf S^{\text{D}} = \mu (\overline{\hat{\mathbf{B}}_{\text{e}}})^{\text{D}}$, and
the evolution equation \eqref{SimMie82} is reduced to
\begin{equation}\label{SimMie8}
-\text{\calligra{L}}_{v}  (\hat{\mathbf{B}}_{\text{e}})
\hat{\mathbf{B}}_{\text{e}}^{-1}
 = \frac{\mu}{\eta} (\overline{\hat{\mathbf{B}}_{\text{e}}})^{\text{D}}.
\end{equation}

\section{Time stepping algorithm}

\subsection{Explicit update formula in Lagrangian formulation}

Let us consider a typical time interval $(t_n, t_{n+1})$ with $\Delta t:= t_{n+1} - t_n > 0$.
By ${}^n \mathbf{C}_{\text{i}}$ and ${}^{n+1} \mathbf{C}_{\text{i}}$ we denote numerical
solutions respectively at $t_n$ and $t_{n+1}$.
Suppose that the deformation gradient ${}^{n+1} \mathbf F$ at time instance $t_{n+1}$
is known, and ${}^n \mathbf{C}_{\text{i}} \in \mathbb{M}$ is given.
The unknown ${}^{n+1} \mathbf{C}_{\text{i}} \in \mathbb{M}$
is estimated as the unimodular part of the solution provided by the classical Euler-backward method (EBM).
In other words, let ${}^{n+1} \mathbf{C}^{\text{EBM}}_{\text{i}} $ be the EBM solution.
After the subsequent correction, the solution is given by
\begin{equation}\label{UniEBM}
{}^{n+1} \mathbf{C}_{\text{i}}  := \overline{{}^{n+1} \mathbf{C}^{\text{EBM}}_{\text{i}}}
= \big(\det({}^{n+1} \mathbf{C}^{\text{EBM}}_{\text{i}})\big)^{-1/3} \ {}^{n+1} \mathbf{C}^{\text{EBM}}_{\text{i}}.
\end{equation}
Let us derive this solution. First, adopting EBM, the discretized version of \eqref{puba44} reads
\begin{equation}\label{DiskrEvol}
{}^{n+1} \mathbf{C}^{\text{EBM}}_{\text{i}}  = {}^{n} \mathbf{C}_{\text{i}} + \frac{\Delta t \mu}{\eta}
\big( {}^{n+1} \overline{\mathbf C}
 ({}^{n+1} \mathbf C^{\text{EBM}}_{\text{i}})^{-1}  \big)^{\text{D}} \ {}^{n+1} \mathbf C^{\text{EBM}}_{\text{i}},
\end{equation}
where ${}^{n+1} \mathbf{C} = {}^{n+1} \mathbf{F}^{\text{T}} \ {}^{n+1} \mathbf{F}$ is given.
Introducing abbreviation
\begin{equation}\label{Abbr}
\beta := \frac{1}{3} \frac{\Delta t \mu}{\eta} \ \text{tr}({}^{n+1} \overline{\mathbf C}
 ({}^{n+1} \mathbf C^{\text{EBM}}_{\text{i}})^{-1}),
\end{equation}
  \eqref{DiskrEvol} can be rewritten as
\begin{equation*}\label{DiskrEvol2}
{}^{n+1} \mathbf{C}^{\text{EBM}}_{\text{i}}  = {}^{n} \mathbf{C}_{\text{i}} + \frac{\Delta t \mu}{\eta}
{}^{n+1} \overline{\mathbf C} - \beta \ {}^{n+1} \mathbf C^{\text{EBM}}_{\text{i}}.
\end{equation*}
Thus, obviously
\begin{equation}\label{DiskrEvol3}
{}^{n+1} \mathbf{C}^{\text{EBM}}_{\text{i}}  =
\frac{1}{1+\beta} \big( {}^{n} \mathbf{C}_{\text{i}} + \frac{\Delta t \mu}{\eta}
{}^{n+1} \overline{\mathbf C} \big).
\end{equation}
Finally, noting that
$\overline{\frac{1}{1+\beta} \mathbf{A}} =
\overline{\mathbf{A}}$ for any $\mathbf{A}$, and substituting \eqref{DiskrEvol3}
into \eqref{UniEBM}, the following explicit update formula is obtained
\begin{equation}\label{ExplUpdate}
\boxed{
{}^{n+1} \mathbf{C}_{\text{i}}  =
\overline{{}^{n} \mathbf{C}_{\text{i}} + \frac{\Delta t \mu}{\eta}
{}^{n+1} \overline{\mathbf C}}}.
\end{equation}
\textbf{Remark 2} Note that the solution ${}^{n+1} \mathbf{C}_{\text{i}}$
is obtained directly, without computing ${}^{n+1} \mathbf{C}^{\text{EBM}}_{\text{i}}$.
On the other hand, a simple explicit update formula can be derived
for EBM as well (see Appendix A).

\subsection{Properties of the algorithm in Lagrangian formulation}

Recall that the exact solution to the evolution equation \eqref{puba44}
has under proper initial conditions
the geometric property $\mathbf{C}_{\text{i}} \in \mathbb{M}$.
The numerical algorithms which exactly preserve this property are referred to as geometric integrators.
It was shown in \cite{ShutovKrStab} that such integrators
allow to suppress the error accumulation, which is especially
important for the simulation of ``long" processes. Therefore, geometric integrators
are advantageous. Obviously, the geometric property \eqref{geopro} is exactly
retained by the numerical solution \eqref{ExplUpdate}.

\textbf{Remark 3} Observe that the method \eqref{UniEBM} corresponds to EBM with a subsequent correction.
Another modification of the classical EBM was considered by Helm  \cite{Helm2} in
order to obtain a geometric integrator. Furthermore, in the paper by Vladimirov et al. \cite{Vladimirov}, two other modifications were considered and
the plastic incompressibility was enforced by introducing the additional equation $\det \mathbf C_{\text{i}} = 1$
with an additional unknown scalar variable. In contrast
to the explicit update formula \eqref{ExplUpdate}, a local iterative
procedure was implemented in \cite{Helm2} and \cite{Vladimirov}.

Another property of the algorithm is as follows.
As already mentioned in Section 2.1, the exact solution  ${\mathbf C}_{\text{i}}$ must be positive definite. Since
${}^{n} \mathbf{C}_{\text{i}}$ and ${}^{n+1} \overline{\mathbf C}$ are positive definite,
the numerical solution ${}^{n+1} \mathbf{C}_{\text{i}}$ given
by \eqref{ExplUpdate} is positive definite as well. Indeed,
the sum of two positive definite tensors is positive definite, and the projection
operator $\overline{(\cdot)}$ retains this property.
Furthermore, due to the positive definiteness, $\det ({}^{n} \mathbf{C}_{\text{i}} + \frac{\Delta t \mu}{\eta}
{}^{n+1} \overline{\mathbf C}) >0$. Thus, the
unimodular part of $({}^{n} \mathbf{C}_{\text{i}} + \frac{\Delta t \mu}{\eta}
{}^{n+1} \overline{\mathbf C})$ in \eqref{ExplUpdate}
is well defined for all $\Delta t \geq 0$ and
${}^{n+1} \mathbf{C}_{\text{i}}$ it is a smooth function of $\Delta t$.
For $\Delta t \geq 0$, the solution ${}^{n+1} \mathbf{C}_{\text{i}}$ ranges
smoothly from ${}^{n} \mathbf{C}_{\text{i}}$ to ${}^{n+1} \overline{\mathbf C}$.
Thus, there is no danger of ``wild oscillations" of the solution,
which are typical for explicit time stepping schemes dealing with large time steps.
The \emph{algorithm is unconditionally stable} since the solution remains
finite for arbitrary time steps.
Like the classical EBM, the algorithm is first order accurate
(see Appendix B).

In case of a relaxation process with $\mathbf C = const$, the Clausius-Duhem inequality \eqref{cld}
requires that the free energy is a decreasing function of time.
Let us analyze the dissipation properties of the presented algorithm in case of stress relaxation.
We consider a local relaxation process with $\mathbf C = const$
and $\mathbf{C}_{\text{i}} = {}^{n} \mathbf{C}_{\text{i}}$ as prescribed initial condition at $t_n$.
In that case, it can be shown that $\psi(\mathbf C \ {{}^{n+1}\mathbf C_{\text{i}}}^{-1})$
is a monotonically decreasing function of $\Delta t$ (see Appendix C).
Thus, the relaxation property of the exact solution is qualitatively reproduced by
the numerical solution.

Note that the evolution equations \eqref{puba44} are invariant
under isochoric change of the reference configuration
(see \cite{ShutovKrStab, ShutovPfeIhl}). More precisely,
let $\mathbf{F}_{0}$ be a constant tensor such that $\det \mathbf{F}_{0} = 1$.
The new reference configuration is introduced as
$\widetilde{K}^{\text{new}}=\mathbf{F}_{0} \widetilde{K}$. Thus, the
corresponding deformation gradient (relative deformation gradient)
is then given by $\mathbf{F}^{\text{new}} = \mathbf{F} \mathbf{F}_0^{-1}$.
Along with the ``old" quantities $\mathbf{C}$ and $\mathbf{C}_{\text{i}}$, we consider
their new counterparts
\begin{equation*}\label{transform}
\mathbf{C}^{\text{new}} = \mathbf{F}_0^{-\text{T}} \mathbf{C} \mathbf{F}_0^{-1}, \quad
\mathbf{C}_{\text{i}}^{\text{new}} =
\mathbf{F}_0^{-\text{T}} \mathbf{C}_{\text{i}} \mathbf{F}_0^{-1}.
\end{equation*}
In that case, it can be easily shown that the evolution of
$\mathbf{C}_{\text{i}}^{\text{new}}$ is governed by
\eqref{puba44} if all ``old" variables are formally replaced
by their new counterparts. The same invariance requirement also
can be formulated for the time-stepping procedure:
For a fixed time step $\Delta t$, we denote the numerical solution by
${}^{n+1}\mathbf C_{\text{i}} = \mathfrak{N}
({}^{n+1}\mathbf C, {}^{n}\mathbf C_{\text{i}})$.
The numerical scheme $\mathfrak{N}$ is called invariant under the reference change if
\begin{equation}\label{transform2}
\mathfrak{N} (\mathbf{F}_0^{-\text{T}} \ {}^{n+1}\mathbf C \ \mathbf{F}_0^{-1},
\mathbf{F}_0^{-\text{T}} \ {}^{n}\mathbf C_{\text{i}} \ \mathbf{F}_0^{-1}) =
\mathbf{F}_0^{-\text{T}} \mathfrak{N} ({}^{n+1}\mathbf C, {}^{n}\mathbf C_{\text{i}}) \mathbf{F}_0^{-1}.
\end{equation}
This property can be checked for the scheme \eqref{ExplUpdate}. Indeed, since $\det \mathbf{F}_{0} = 1$, we have
\begin{equation*}\label{transform3}
\overline{\mathbf{F}_0^{-\text{T}} {}^{n} \mathbf{C}_{\text{i}} \mathbf{F}_0^{-1} + \frac{\Delta t \mu}{\eta} \
\overline{\mathbf{F}_0^{-\text{T}} \ {}^{n+1}\mathbf C \ \mathbf{F}_0^{-1}}} =
\mathbf{F}_0^{-\text{T}} \
\overline{{}^{n} \mathbf{C}_{\text{i}} + \frac{\Delta t \mu}{\eta} \
{}^{n+1} \overline{\mathbf C}} \
\mathbf{F}_0^{-1}.
\end{equation*}

\textbf{Remark 4} One possible application of the invariance property \eqref{transform2}
is as follows. By choosing $\mathbf{F}_{0} := ({}^{n} \mathbf{C}_{\text{i}})^{1/2}$ we
get ${}^{n} \mathbf{C}^{\text{new}}_{\text{i}}=\mathbf{1}$. Thus, any algorithm computing $\mathfrak{N} ({}^{n+1}\mathbf C, \mathbf{1})$
would be sufficient to restore the general scheme $\mathfrak{N}$. In other words, it is
sufficient to perform time-stepping for zero initial conditions for the inelastic strain, which corresponds
to ${}^{n} \mathbf{C}_{\text{i}} = \mathbf{1}$.

Within implicit, deformation driven finite element procedures, the
consistent tangent operator is required for iterative solving
for global balance of linear momentum \cite{SimHug}.
An explicit expression for the consistent tangent is presented in Appendix D, and
its symmetry is proved.

It can be shown that the numerical solution \eqref{ExplUpdate}
corresponds to the \emph{exact solution} of a relaxation problem with a
fixed $\mathbf{C}$ at a certain time instance $t \approx t_{n+1}$.
More precisely, consider the initial
value problem \eqref{puba44}, \eqref{initCond} with $\mathbf{C}(t) = const$.
Then the exact solution can be represented as
$\mathbf C_{\text{i}}(t) = \overline{\mathbf C_{\text{i}}^0 + \varphi(t)
\overline{\mathbf C}}$ with $\varphi(t^0) = 0$ and
$\varphi(t) \approx \frac{\mu (t-t^0)}{\eta}$
as $t \rightarrow t^0$ (see Appendix E).

\subsection{Explicit update formula in Eulerian formulation}

In this subsection we rewrite the explicit
update formula \eqref{ExplUpdate} in terms of tensors which operate on
the current configuration. Toward that end, for each time interval $(t_n, t_{n+1})$, we introduce
the so-called trial elastic left Cauchy-Green tensor under the assumption of a frozen inelastic flow.
Having in mind that ${\hat{\mathbf B}}_{\text{e}} = \mathbf F \ \mathbf{C}^{-1}_{\text{i}} \ {\mathbf F}^{\text{T}}$,
we get the trial elastic strain in the form
\begin{equation*}\label{TrialBe}
{}^{n+1} {\hat{\mathbf B}}^{\text{trial}}_{\text{e}} =
{}^{n+1} {\mathbf F} \  {}^{n} \mathbf{C}^{-1}_{\text{i}} \ {}^{n+1} {\mathbf F}^{\text{T}}.
\end{equation*}
Multiplying \eqref{ExplUpdate} with ${}^{n+1} \overline{{\mathbf F}}^{-\text{T}}$ from the left and with
${}^{n+1} \overline{{\mathbf F}}^{-1}$ from the right, we get
\begin{equation*}\label{EulerStep0}
{}^{n+1} \overline{{\hat{\mathbf B}}}^{-1}_{\text{e}} =
\overline{({}^{n+1}\overline{{\hat{\mathbf B}}_{\text{e}}^{\text{trial}}})^{-1} + \frac{\Delta t \mu}{\eta} \ \mathbf{1}}.
\end{equation*}
Taking the inelastic incompressibility into account, one gets
\begin{equation}\label{EulerStep}
\boxed{
{}^{n+1} {\hat{\mathbf B}}^{-1}_{\text{e}} = (\det {}^{n+1} \mathbf F)^{-2/3}
\overline{({}^{n+1}\overline{{\hat{\mathbf B}}_{\text{e}}^{\text{trial}}})^{-1} +\frac{\Delta t \mu}{\eta} \ \mathbf{1}}}.
\end{equation}
This is an \emph{explicit} update formula for the evolution equation \eqref{SimMie8}. Note that,
similar to the product formula consistent with the operator split considered by Simo  \cite{Simo},
${}^{n+1} {\hat{\mathbf B}}_{\text{e}}$ is co-axial with ${}^{n+1} {\hat{\mathbf B}}^{\text{trial}}_{\text{e}}$.
Within a time step, the explicit update formula \eqref{EulerStep} predicts
the same stress response as the formula \eqref{ExplUpdate}.

\subsection{Exponential mapping: Lagrangian and Eulerian formulations}

Alternatively to the explicit update formula considered in this paper,
a well-known implicit method based on the exponential mapping can be adopted:
\begin{equation}\label{ExpMeth}
{}^{n+1} {\mathbf C}_{\text{i}} =
\exp\big[\frac{\Delta t \mu}{\eta} ({}^{n+1} \overline{\mathbf C} \ {}^{n+1}
\mathbf C^{-1}_{\text{i}})^{\text{D}}\big] {}^{n} {\mathbf C}_{\text{i}}.
\end{equation}
By some algebraic computations, it can be shown that the invariance relation  \eqref{transform2} holds
for this method as well. Next, a proof that the symmetry property is retained by the algorithm
even in case of a more general material model is presented in \cite{ShutovKrVisc}.
Moreover, the incompressibility of the inelastic flow is exactly retained.
Thus, for the corresponding numerical solution, ${}^{n+1} {\mathbf C}_{\text{i}} \in \mathbb{M}$.

Observe that the numerical scheme obtained
using the operator split technique (cf. Simo  \cite{Simo})
can be derived directly from \eqref{ExpMeth}.
Indeed, by inverting both sides of
\eqref{ExpMeth}, we get
\begin{equation}\label{ExpMeth2}
{}^{n+1} {\mathbf C}^{-1}_{\text{i}} = {}^{n} {\mathbf C}^{-1}_{\text{i}} \
\exp\big[- \frac{\Delta t \mu}{\eta} ({}^{n+1} \overline{\mathbf C} \ {}^{n+1}
\mathbf C^{-1}_{\text{i}})^{\text{D}}\big] .
\end{equation}
Since the inelastic incompressibility is exactly retained by \eqref{ExpMeth},
$\det ({}^{n+1} {\hat{\mathbf B}}^{-1}_{\text{e}}) = (\det {}^{n+1} \mathbf F)^2$. Therefore
\begin{equation}\label{ExpMeth3}
{}^{n+1} {\mathbf F}^{-\text{T}} \ ({}^{n+1} \overline{\mathbf C} \ {}^{n+1}
\mathbf C^{-1}_{\text{i}})^{\text{D}} \ {}^{n+1} {\mathbf F}^{-\text{T}} =
({}^{n+1} \overline{\hat{\mathbf B}}_{\text{e}})^{\text{D}}.
\end{equation}
Multiplying \eqref{ExpMeth2} with ${}^{n+1} {\mathbf F}$ from the left
and with ${}^{n+1} {\mathbf F}^{\text{T}}$ from the right, we get using \eqref{ExpMeth3}
\begin{equation*}\label{ExpMeth4}
{}^{n+1} {\hat{\mathbf B}}_{\text{e}} = {}^{n+1} {\hat{\mathbf B}}^{\text{trial}}_{\text{e}} \
\exp\big[- \frac{\Delta t \mu}{\eta} ({}^{n+1} \overline{\hat{\mathbf B}}_{\text{e}})^{\text{D}}\big].
\end{equation*}
Therefore, ${}^{n+1} {\hat{\mathbf B}}_{\text{e}} $ commutes with ${}^{n+1} {\hat{\mathbf B}}^{\text{trial}}_{\text{e}}$.
Thus, this equation can be rewritten in the form
\begin{equation}\label{ExpMeth5}
{}^{n+1} {\hat{\mathbf B}}_{\text{e}} =
\exp\big[- \frac{\Delta t \mu}{\eta}
({}^{n+1} \overline{\hat{\mathbf B}}_{\text{e}})^{\text{D}}\big] \ {}^{n+1} {\hat{\mathbf B}}^{\text{trial}}_{\text{e}} .
\end{equation}
The well-known implicit update formula is restored (see, for instance, equation (44) in \cite{ReeseG}).
We stress that the nonlinear equations \eqref{ExpMeth} and \eqref{ExpMeth5}
represent one and the same method, written in two different formulations.
For a given time step, these equations are equivalent since they predict
the same stress response.

\section{Numerical results}

\subsection{Accuracy testing for a single Maxwell element}

\begin{figure}\centering
\psfrag{A}[m][][1][0]{$\| {\mathbf C}^{numer}_{\text{i}} - {\mathbf C}^{exact}_{\text{i}} \|$}
\psfrag{B}[m][][1][0]{$\Delta t = 0.5$}
\psfrag{C}[m][][1][0]{$\Delta t = 1$}
\scalebox{0.8}{\includegraphics{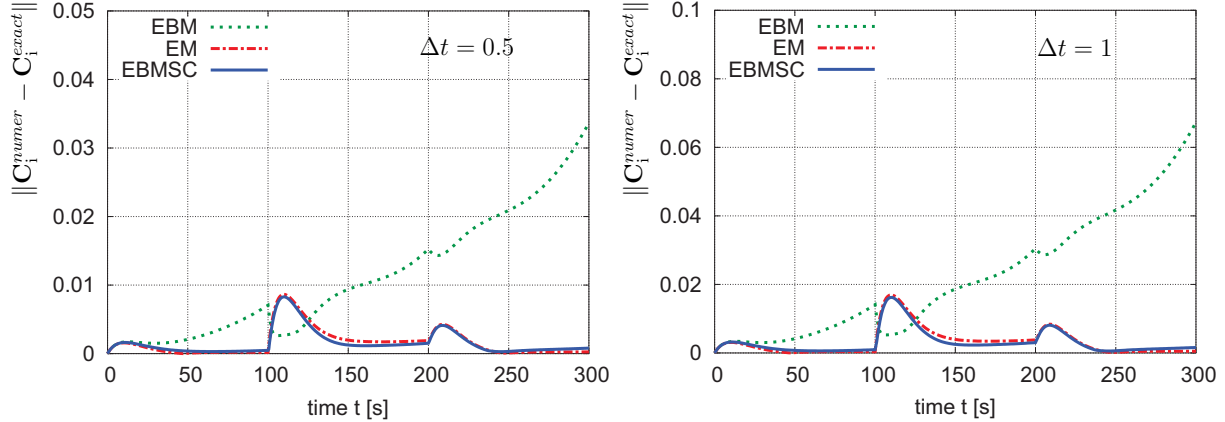}}
\caption{Plots of the numerical error pertaining to three different
algorithms:
classical Euler backward (EBM),
exponential method (EM), and the
Euler backward with subsequent correction (EBMSC).
Explicit update formula for EBMSC is given by \eqref{ExplUpdate}. All methods are first order accurate.
Geometric integrators (EBMSC and EM) prevent the error accumulation.  \label{fig2}}
\end{figure}

\begin{figure}\centering
\psfrag{A}[m][][1][0]{$\mathbf{T}_{1 1}$ [MPa]}
\psfrag{B}[m][][1][0]{$\mathbf{T}_{1 2}$ [MPa]}
\psfrag{D}[m][][1][0]{$\Delta t = 10$}
\psfrag{E}[m][][1][0]{$\Delta t = 5$}
\psfrag{F}[m][][1][0]{$\Delta t = 1$}
\scalebox{0.8}{\includegraphics{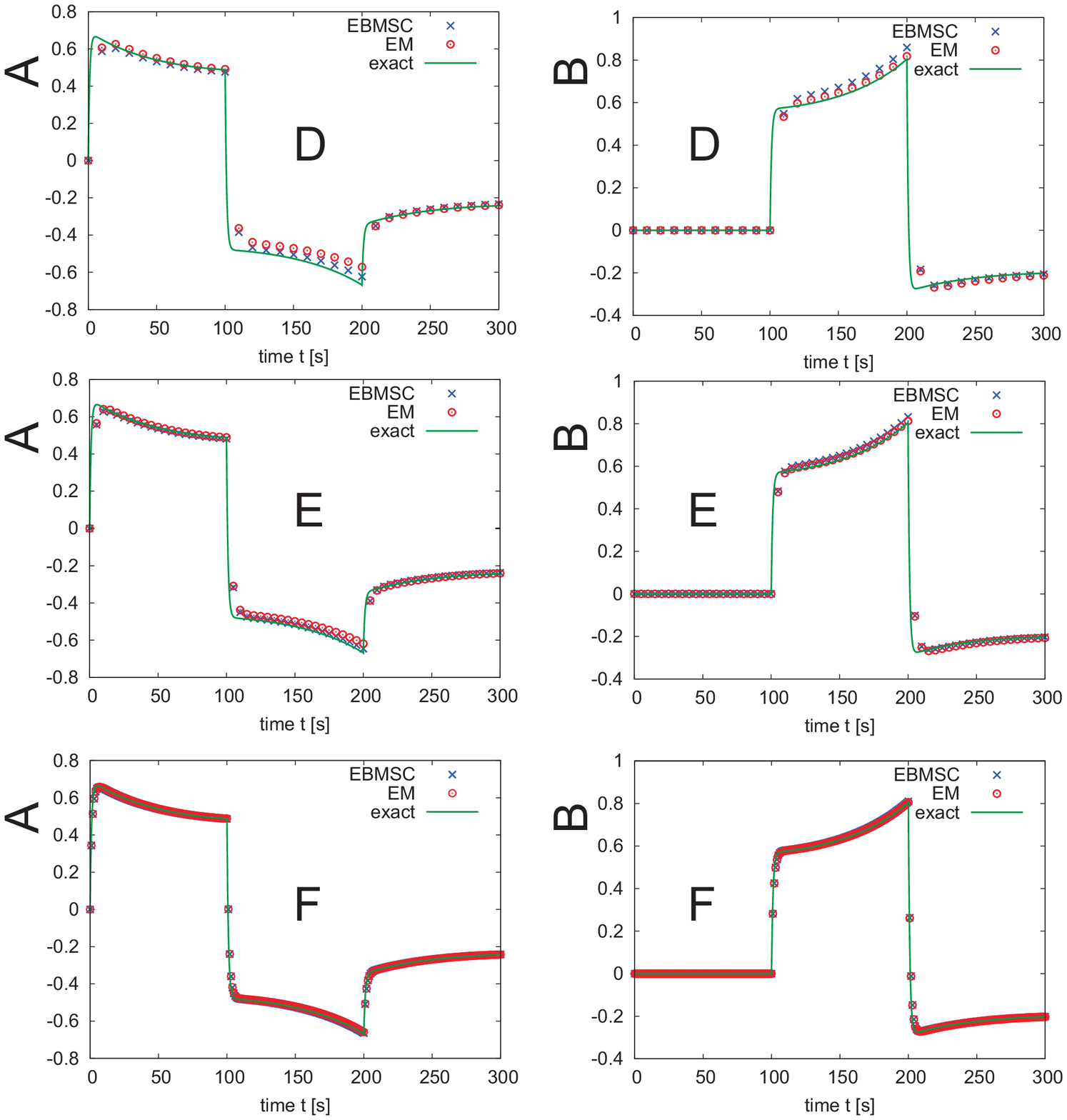}}
\caption{Simulated stress response in terms of Cauchy stresses for $\Delta t = 10$ (top),
$\Delta t = 5$ (middle), and $\Delta t = 1$ (bottom).
Since both methods are first order accurate,
the error is proportional to $\Delta t$.
  \label{fig3}}
\end{figure}

In order to test the accuracy of the
explicit update formula \eqref{ExplUpdate}, we consider a
local loading program in the time interval $t \in [0,300]$ (time is measured in seconds)
\begin{equation}\label{loaprog0}
\mathbf F (t) = \overline{\mathbf F^{\prime} (t)},
\end{equation}
where $\mathbf F^{\prime}(t)$ is a piecewise linear function of time $t$ such that
$\mathbf F^{\prime} (0) = \mathbf F_1$, $\mathbf F^{\prime} (100) = \mathbf F_2$,
$\mathbf F^{\prime} (200) = \mathbf F_3$, and $\mathbf F^{\prime} (300) = \mathbf F_4$ with
\begin{equation*}\label{loaprog2}
\mathbf F_1 :=\mathbf 1, \quad
\mathbf F_2 := 2 \mathbf{e}_{1} \otimes \mathbf{e}_{1} +
\frac{1}{\sqrt2} (\mathbf{e}_{2} \otimes \mathbf{e}_{2} + \mathbf{e}_{3} \otimes \mathbf{e}_{3}),
\end{equation*}
\begin{equation*}\label{loaprog222}
\mathbf F_3 := \mathbf 1 + \mathbf{e}_{1} \otimes \mathbf{e}_{2}, \quad
\mathbf F_4 := 2 \mathbf{e}_{2} \otimes \mathbf{e}_{2} +
\frac{1}{\sqrt2} (\mathbf{e}_{1} \otimes \mathbf{e}_{1} + \mathbf{e}_{3} \otimes \mathbf{e}_{3}).
\end{equation*}
More precisely, we put here
\begin{equation*}\label{loaprog}
\mathbf F^{\prime} (t) :=
\begin{cases}
    (1 - t/100) \mathbf F_1  + (t/100) \mathbf F_2 \quad \quad  \ \  \text{if} \ t \in [0,100] \\
    (2 - t/100) \mathbf F_2  + (t/100-1) \mathbf F_3 \quad \text{if} \ t \in (100,200] \\
     (3 - t/100) \mathbf F_3  + (t/100-2) \mathbf F_4 \quad \text{if} \ t \in (200,300]
\end{cases}.
\end{equation*}
The reference configuration is assumed to be stress free at $t=0$. Thus, we put
$\mathbf C_{\text{i}}|_{t=0} =  \mathbf 1$.
The following values of the
material parameters are used: $\eta = 400 \ \text{MPa s}$, $\mu = 40$ MPa.
The numerical solution of the initial value problem
obtained with extremely small time step
($\Delta t = 0.001$) will be referred to as
the exact solution ${\mathbf C}^{exact}_{\text{i}}$.
The numerical solutions with
$\Delta t = 1$ and $\Delta t = 0.5$ are
denoted by ${\mathbf C}^{numer}_{\text{i}}$. The error
$\| {\mathbf C}^{numer}_{\text{i}} - {\mathbf C}^{exact}_{\text{i}} \|$
is plotted versus time in Fig. \ref{fig2} for three different methods.
The explicit update \eqref{ExplUpdate} is abbreviated as
Euler-backward method with subsequent correction (EBMSC).
Additionally to EBMSC, the classical EBM method and the exponential method (EM)
are considered in this subsection.
Since these three methods are first order accurate, the error is proportional to $\Delta t$.
Moreover, in accordance with \cite{ShutovKrStab}, there is no error accumulation
in case of geometric integrators (EBMSC and EM). More precisely,
the error is uniformly bounded by $C \Delta t$ where the constant $C$ does not depend on the
size of the entire time interval \cite{ShutovKrStab}.
Next, since the incompressibility condition is
violated by EBM, the geometric property \eqref{geopro}
is not preserved and the numerical error tends to accumulate over time.
The geometric integrators EBMSC and EM are equivalent in terms of accuracy.
The simulated stress response is presented in  Fig. \ref{fig3} for three different time steps $\Delta t$.

\begin{figure}\centering
\psfrag{A}[m][][1][0]{$\mathbf{T}_{1 1}$ [MPa]}
\psfrag{B}[m][][1][0]{$\mathbf{T}_{1 2}$ [MPa]}
\psfrag{D}[m][][1][0]{$\eta = 4 \cdot 10^{2} \text{MPa s}$}
\psfrag{E}[m][][1][0]{$\eta = 4 \cdot 10^{3} \text{MPa s}$}
\psfrag{F}[m][][1][0]{$\eta = 4 \cdot 10^{4} \text{MPa s}$}
\psfrag{K}[m][][1][0]{a)}
\psfrag{L}[m][][1][0]{b)}
\psfrag{M}[m][][1][0]{c)}
\scalebox{0.75}{\includegraphics{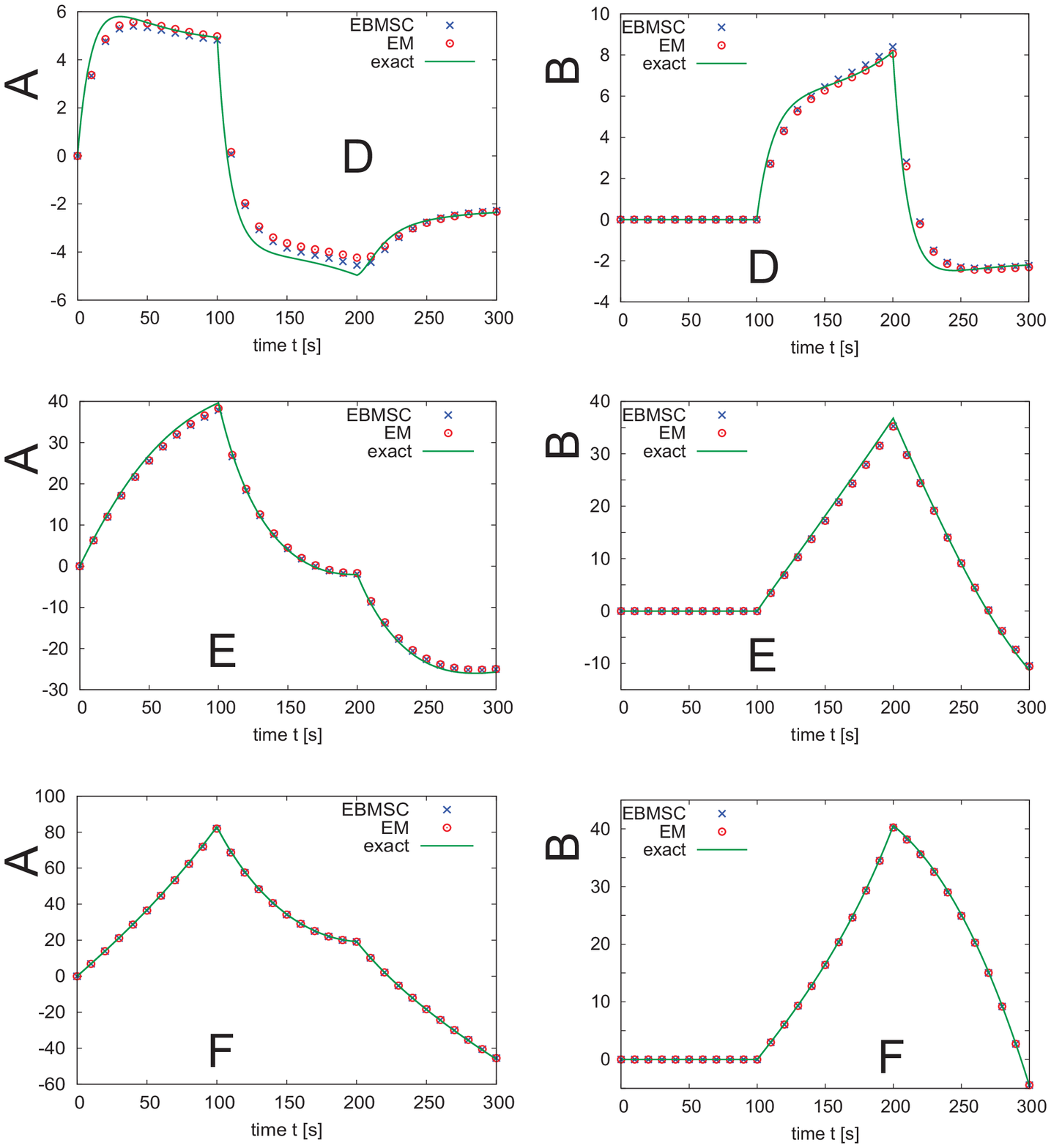}}
\caption{Simulated stress response in terms of Cauchy stresses for
different values of $\eta$: $\eta = 4 \cdot 10^{2} \text{MPa s}$ (top),
$\eta = 4 \cdot 10^{3} \text{MPa s}$ (middle), and
$\eta = 4 \cdot 10^{4} \text{MPa s}$ (bottom).
The relative error is higher for stiffer problems with small $\eta$.
  \label{fig4}}
\end{figure}

Next, the accuracy is tested for varying viscosity
 $\eta$ with a fixed $\mu = 40$ MPa, $\Delta t = 10$
and loading programm given by \eqref{loaprog0}.
The largest relative error is observed for the smallest
values of $\eta$, although the overall stress level tends to zero (see the top part of Fig. \ref{fig4}).
Indeed, for smaller $\eta$ the problems becomes stiffer. Although the implicit schemes
are robust, smaller time steps are still required to reduce the
numerical error. In case of large $\eta$, the inelastic flow is
insignificant and the material response becomes
nearly hyperelastic (see the bottom part of Fig. \ref{fig4}).

\subsection{Application to visoelasticity within the FEM}

\subsubsection{Material model of finite strain viscoelasticity}

%A viscoelastic model of a rubber-like material
%is considered in this subsection.
In this section, a practical application of the proposed
integration algorithm to the modeling of viscoelastic material response
is discussed. As a typical example,
a constitutive model proposed by Reese and Govindjee \cite{ReeseG, GovRes}
is considered. The rheological interpretation of the model consists of
a single hyperelastic spring (Hooke body) for the representation of equilibrium stresses
and $N$ Maxwell bodies which are connected in parallel to the Hooke body to represent
viscous effects, see Fig. \ref{fig1}b.
Each of these Maxwell bodies is governed by equations presented in Section 2.
The total free energy is given as a sum of isotropic functions as follows (cf. \cite{ReeseG, GovRes}):
\begin{equation*}\label{freeEnergyRubber}
    \psi =  \psi_{\text{eq}}(\mathbf{C}) + \sum\limits_{m=1}^{N} \psi_{\text{ov},m}(\mathbf{C} \mathbf{C}_{\text{i},m}^{-1}),
\end{equation*}
where the equilibrium part is modeled by
the hyperelastic Yeoh model \cite{Yeoh 1993} in combination with an additional
volumetric part
\begin{equation}\label{Equilibrium}
 \rho_{\scriptscriptstyle \text{R}} \psi_{\text{eq}}(\mathbf{C}) =
 \sum_{n=1}^3  c_{i0} \Big({\rm tr} \overline{\mathbf{C}} - 3\Big)^n
 + \frac{9\,k}{2} \Big( (\det \mathbf{F})^{1/3} - 1\Big)^2.
\end{equation}
Here, $c_{10}$, $c_{20}$, and $c_{30}$ are material parameters of the Yeoh material;
 $k$ stands for the bulk modulus.
Observe that the volumetric stress response depends
on the volume ratio, given by  $\det\mathbf{F}$.\footnote{
The volumetric ansatz presented in \eqref{Equilibrium} is implemented in MSC.MARC \cite{MARC}.
Obviously, any alternative ansatz for the volumetric part can be used as well.}
Next, the free energy for the $m^{th}$ Maxwell body reads as
\begin{equation*}\label{Noneq}
\rho_{\scriptscriptstyle \text{R}} \psi_{\text{ov},m}=
\rho_{\scriptscriptstyle \text{R}} \psi_{\text{ov},m}(\mathbf C \mathbf C_{\text{i},k}^{-1}) =
\frac{\mu_m}{2 } \big(
\text{tr} \overline{\mathbf C \mathbf C_{\text{i},m}^{-1}} - 3 \big), \quad m=1,2,...,N.
\end{equation*}
Here, $\mu_m \geq 0$ is the shear modulus of the $m^{th}$ element and $\mathbf{C}_{\text{i},m}$
is the corresponding inelastic tensor of right Cauchy Green type.
The evolution of each of these variables is governed by equations of type \eqref{puba44}.
The resulting system of constitutive equations is as follows:
\begin{equation*}\label{SumStresses}
  \tilde{\mathbf{T}}
      = \tilde{\mathbf{T}}_{\text{eq}} + \sum\limits_{m=1}^{N} \tilde{\mathbf{T}}_{\text{ov},m},
\end{equation*}
\begin{equation*}\label{EquilibriumStress}
  \tilde{\mathbf{T}}_{\text{eq}}
  = \Big(2\ c_{10} + 4\,c_{20}\big({\rm tr}\overline{\mathbf{C}}-3\big) +
  6 \ c_{30}\big({\rm tr}\overline{\mathbf{C}}-3\big)^2\Big)\ \overline{\mathbf{C}}^{\text{D}}\mathbf{C}^{-1}
  + 3 k (\det \mathbf{F})^{1/3} \ \big((\det \mathbf{F})^{1/3} - 1\big) \ \mathbf{C}^{-1},
\end{equation*}
\begin{equation*}\label{eq:constitutive_equations_full_model}
  \tilde{\mathbf{T}}_{\text{ov},m}
  = \mu_m \mathbf{C}^{-1} \big(\overline{\mathbf{C}}\mathbf{C}_{\text{i},m}^{-1}\big)^{\text{D}},  \quad
  \dot{\mathbf{C}}_{\text{i},m}
  =  \frac{\mu_m}{\eta_m} \big(\overline{\mathbf{C}}\mathbf{C}_{\text{i},m}^{-1}\big)^{\text{D}}
   \mathbf{C}_{\text{i},m}, \quad m=1,2,...,N,
\end{equation*}
where $\eta_m \geq 0$ denotes the viscosity of the $m^{th}$ Maxwell body.

The material model is implemented into commercial finite element
software MSC.MARC\textsuperscript{\textregistered}, making use of the
user subroutine \textit{HYPELA2}.
Since the total strain is prescribed at each point of
Gauss integration, the numerical treatment of each of the Maxwell elements
is fully independent from other Maxwell elements.
For each of them, the corresponding
evolution equations are solved with the help of the explicit
update formula \eqref{ExplUpdate} and the consistent tangent operator
is computed explicitly (cf. Appendix D).
In order to avoid volume locking effects, a
mixed u-p-formulation is adopted \cite{SussmanBathe1987}.
For the subsequent simulations, a system of 4 Maxwell bodies is implemented with
the material parameters from Tab. \ref{tab1}.
\begin{table}[h]
 \caption{Material parameters of a viscoelastic material ($m=1,2,3,4$).}
	\centering
		\begin{tabular}{|c | c |c |c |c |c |c |c |c |c|}
      \hline
         $c_{10} \, {\rm [MPa]}$
      & $c_{20} \, {\rm [MPa]}$
      & $c_{30} \, {\rm [MPa]}$
      & $k \, {\rm [MPa]}$
      & $\ N \ {\rm [-]}$
      & $\mu_m \, {\rm [MPa]}$
      & $\eta_m \, {\rm [MPa \ s]}$\\
      \hline &&&&&& \\[-8mm]
      0.45 & -0.048 & 0.011 & 1000 & $4$& 0.2 & $2 \cdot 10^{m-3}$ \\
      \hline
		\end{tabular}
    \label{tab1}
\end{table}

\begin{figure}\centering
\psfrag{A}[m][][1][0]{$\mathbf{T}_{x x}$ [MPa]}
\psfrag{B}[m][][1][0]{$\mathbf{F}_{x x} [-]$}
\scalebox{0.8}{\includegraphics{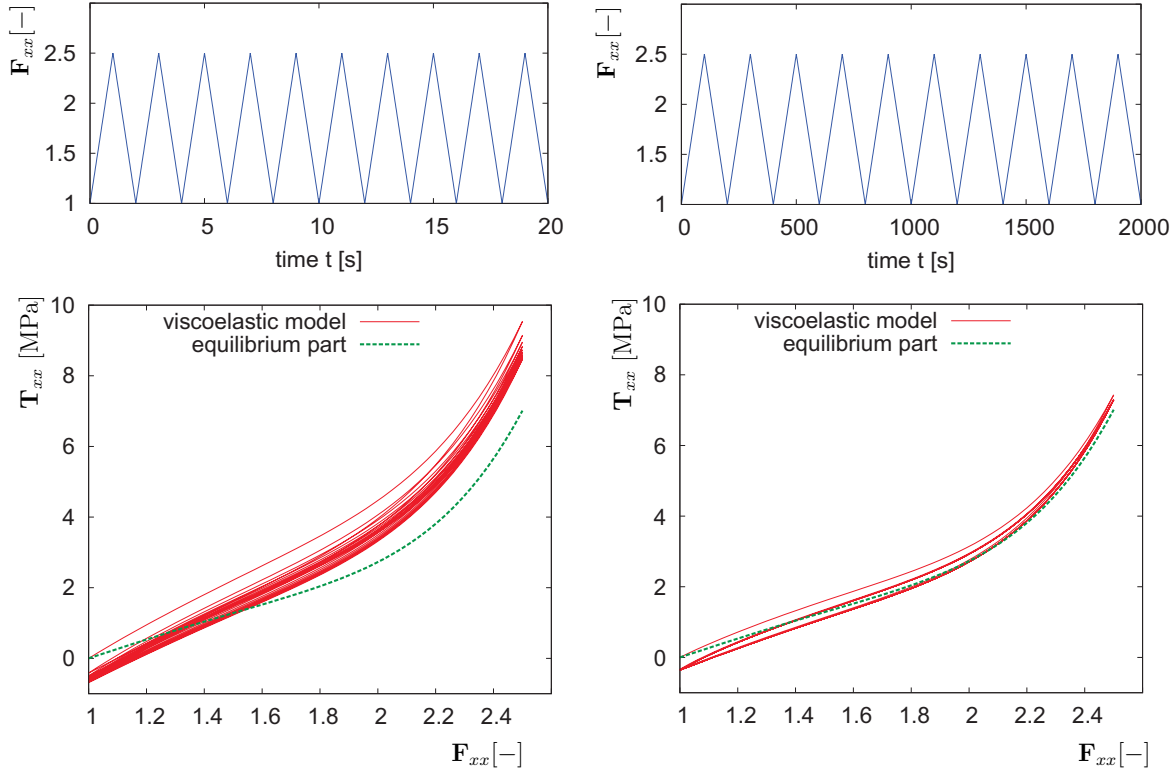}}
\caption{Simulated stress response under uniaxial cyclic loading for
different strain rates: for $ | \dot{\mathbf{F}}_{xx} | = 1.5 \ \text{s}^{-1}$  (left)
and for $| \dot{\mathbf{F}}_{xx} | = 0.015 \ \text{s}^{-1}$ (right).
  \label{fig5}}
\end{figure}

\begin{figure}\centering
\psfrag{A}[m][][1][0]{$\mathbf{T}_{1 1}$ [MPa]}
\psfrag{B}[m][][1][0]{$\mathbf{F}_{1 1} [-]$}
\scalebox{0.8}{\includegraphics{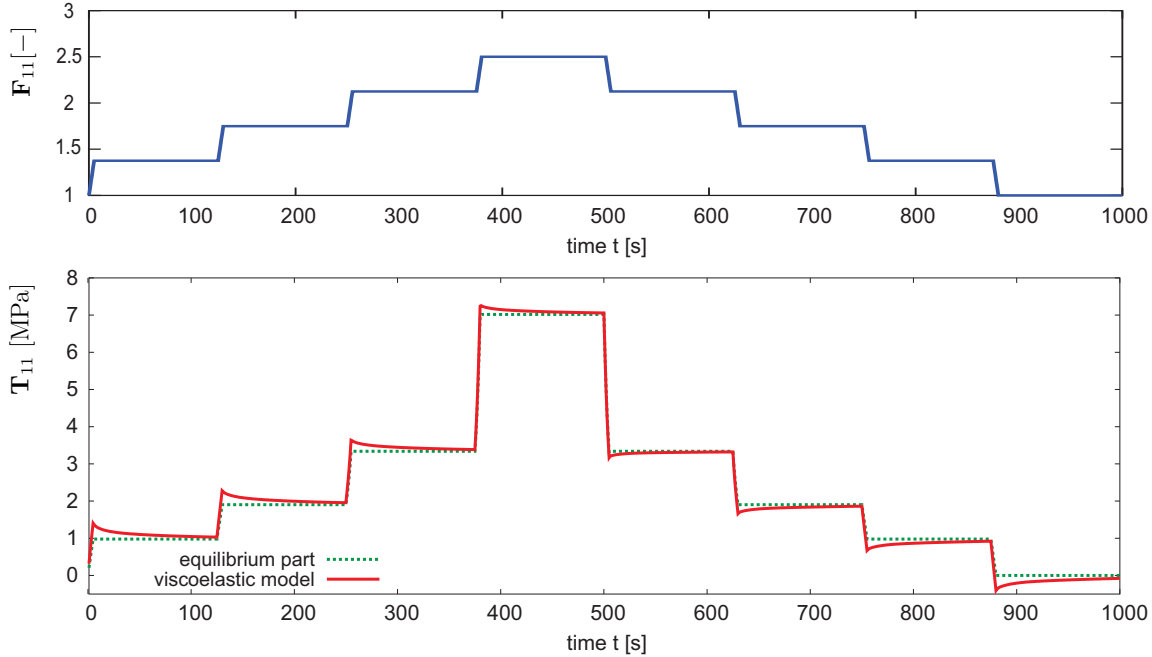}}
\caption{Simulation results for uniaxial relaxation tests: loading program (top) and stress response (bottom).
  \label{fig6}}
\end{figure}

First, in order to illustrate the stress response of the model, a
uniaxial cyclic loading is simulated.
The loading axis is fixed and coincides with the $x$-axis.
In Fig. \ref{fig5} the stress response under a strain driven loading
is shown for two different loading rates.
It can be seen that for small loading rates the stress response converges
to the equilibrium path.
Next, the simulated stress response for uniaxial relaxation test is represented in Fig. \ref{fig6}.
Note that if the instant stress lies under the equilibrium stress, the relaxation process results in
increasing stresses.

\subsubsection{FEM solution of a representative boundary value problem}

\begin{figure}\centering
\includegraphics[width = .8\textwidth]{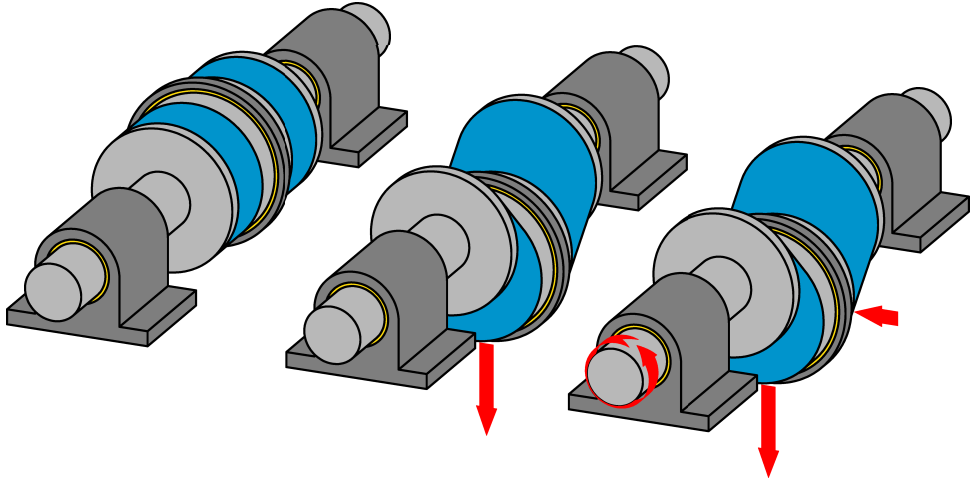}
\caption{Experimental setup for simple shear loading with rotating axes \cite{Klauke1}:
initial state  (left), monotonic shear loading (middle),
shear loading with rotating axes (right).
  \label{fig7}}
\end{figure}

A representative boundary value problem is solved in this subsection in order to
test the stability properties of the proposed integration algorithm.
Toward that end, the material model from the previous subsection is adopted to simulate
the so-called rotary dynamics experiment used for the analysis of
rubber materials. The corresponding experimental setup was originally proposed in
\cite{Gent 1960}
and applied later for the life-time prediction of rubber materials in \cite{Klauke1, Klauke2}.
Within this experiment, a shear loading with rotating axes is
applied to two cylindrical specimens as
shown in Fig. \ref{fig7}.\footnote{The animated version of the experiment can be
seen at http://youtu.be/eNgjGE7upYY .}

\begin{figure}\centering
\scalebox{0.8}{\includegraphics{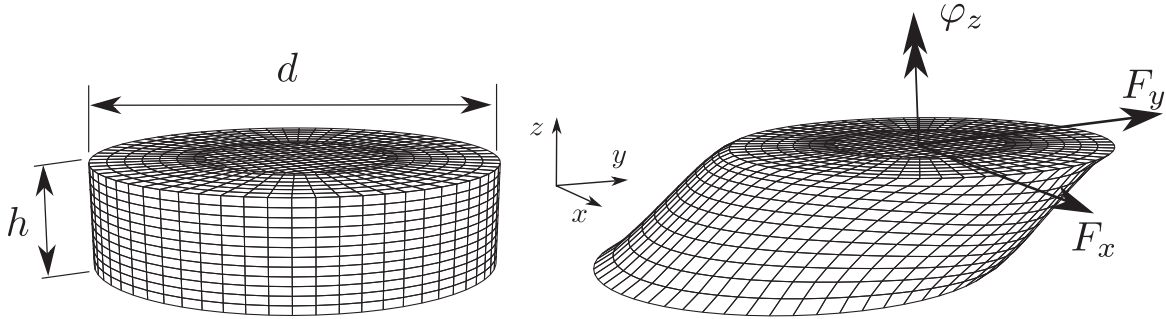}}
\caption{FEM model of the rotary dynamics testing: initial configuration (left) and
deformed configuration (right).
  \label{fig8}}
\end{figure}

The used finite element model
is adopted from \cite{Hohl 2007}. It represents a single cylinder of height
$h = 40$ mm and of diameter $d = 10$ mm (Fig. \ref{fig8}).
Eight-node isoparametric three-dimensional brick elements with trilinear interpolation
and one extra node with a single degree of freedom for pressure (MSC.MARC element type 84) are used.
The mesh consists of 7000 elements. Loads are applied only on the upper and lower surfaces in the following manner.
Both surfaces are rigid and can rotate independently about corresponding rotation axes.
Each rotation axis goes through the center of the corresponding surface in the normal direction.
For the lower surface, the rotation axis is fixed and the rotation is free.
For the upper surface, an initial displacement of the surface is submitted in $y$-direction
as shown in Fig. \ref{fig8},
such that a simple shear loading is applied to the specimen.
In the next loading step, the corresponding axis is fixed and a
rotation of the upper surface is prescribed with a constant
angular velocity  $\dot{\varphi}_z = 0.2 \pi /\text{s}$.
Thus, the sample is subjected to a non-proportional cyclic loading.
Finally, after $300$ s, the rotation is stopped.

\begin{figure}\centering
\scalebox{0.8}{\includegraphics{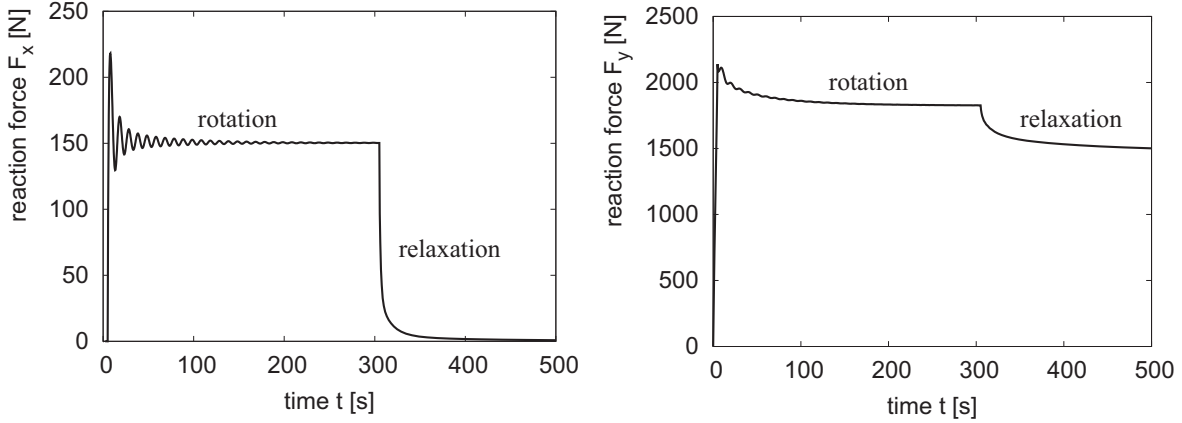}}
\caption{Simulation results for reaction forces.
FEM simulation performed in MSC.MARC using
explicit update formula \eqref{ExplUpdate}.
  \label{fig9}}
\end{figure}

The numerical simulation was performed with the time step $\Delta t = 0.25$ s. The
material parameters of the viscoelastic material are taken from Tab. \ref{tab1}.
No convergence difficulties were observed during the simulation whatsoever.
The simulated reaction forces $F_x$ and $F_y$ are shown in
Fig. \ref{fig9}. After approximately 10 revolutions of the sample, nearly stationary solution is observed.
Finally, during the relaxation step, the reaction force $F_x$ disappears and
the force $F_y$ tends to a certain equilibrium value.

\section{Discussion and conclusion}

A \emph{closed from solution} for fully implicit integration algorithm is proposed.
The algorithm corresponds to the classical Euler-backward method with a subsequent correction
to enforce the incompressibility of the inelastic flow.
The scheme is highly efficient, since no local iteration procedure is required.
The numerical solution is well defined even for large time steps $\Delta t$.
The resulting integration algorithm is unconditionally stable and first order accurate.
The algorithm shows a similar accuracy as the well-known exponential scheme (EM).
The consistent tangent operator is symmetric.
The following properties of the exact solution are retained by the algorithm \eqref{ExplUpdate}:
\begin{itemize}
\item[i] symmetry and incompressibility: ${}^{n+1} \mathbf{C}_{\text{i}} \in Sym$, $\det({}^{n+1} \mathbf{C}_{\text{i}}) = 1$,
\item[ii] positive definiteness: ${}^{n+1}\mathbf{C}_{\text{i}} > 0$,
\item[iii] the solution ${}^{n+1}\mathbf{C}_{\text{i}}$ is a
smooth function of $\Delta t$,
\item[iv] for the stress relaxation with $\mathbf C = const$,
the free energy $\psi(\mathbf C \ {{}^{n+1}\mathbf C_{\text{i}}}^{-1})$
is a monotonically decreasing function of $\Delta t$,
\item[v] numerical solution remains invariant under
the isochoric change of the reference configuration.
\end{itemize}
No error accumulation is observed due to the exact
preservation of the incompressibility constraint, in accordance
to the theoretical results from \cite{ShutovKrStab}.

The explicit update formula is derived for a special variant of the finite strain
Maxwell fluid which is widely adopted in material modelling.
Some modifications of the algorithm are possible
to cover more general material behavior.
In particular, the case of process-dependent viscosity
can be considered. Moreover, the Neo-Hookean type of hyperelasticity
can be replaced by more general constitutive assumptions.
Application of the method to elasto-plasticity with
different types of nonlinear hardening seems promising.

\section*{Acknowledgement}

This research was supported by the German National Science Foundation (DFG) within SFB 692, PAK 273 and SFB/TR 39.

\section*{Appendix A (explicit update formula for EBM)}
Starting from \eqref{DiskrEvol3}, we derive an explicit update
formula for the classical Euler-backward method.
It is sufficient to obtain a
closed-form relation for $\beta$ which appears in \eqref{DiskrEvol3}. First, we abbreviate:
\begin{equation*}\label{AppendixA1}
\mathbf{\Phi} := {}^{n} \mathbf{C}_{\text{i}} +  \frac{\mu \Delta t}{\eta}
{}^{n+1} \overline{\mathbf C}.
\end{equation*}
By computing the inverse of both sides of \eqref{DiskrEvol3}, we get
\begin{equation*}\label{AppendixA2}
( {}^{n+1} \mathbf{C}^{\text{EBM}}_{\text{i}} )^{-1} = (1 + \beta) \mathbf{\Phi}^{-1}.
\end{equation*}
Substituting this result into the definition \eqref{Abbr}, we obtain a linear equation with respect to $\beta$
\begin{equation*}\label{AppendixA3}
\beta = \frac{1}{3} \frac{\mu \Delta t}{\eta}  \text{tr}\big( (1+\beta) \ {}^{n+1} \mathbf{C} \mathbf{\Phi}^{-1}\big).
\end{equation*}
After resolving it, we get
\begin{equation*}\label{AppendixA4}
\beta = \frac{\frac{1}{3} \frac{\mu \Delta t}{\eta} \ \text{tr}({}^{n+1}
\mathbf{C} \mathbf{\Phi}^{-1})}{1-\frac{1}{3} \frac{\mu \Delta t}{\eta}
 \ \text{tr}({}^{n+1} \mathbf{C} \mathbf{\Phi}^{-1})} .
\end{equation*}
Finally, substituting this into \eqref{DiskrEvol3}, the EBM solution is given by
\begin{equation*}\label{AppendixA5}
{}^{n+1} \mathbf{C}^{\text{EBM}}_{\text{i}}  =
\big(1-\frac{1}{3} \frac{\mu \Delta t}{\eta}
\ \text{tr}({}^{n+1} \mathbf{C} \mathbf{\Phi}^{-1})\big) \mathbf{\Phi}.
\end{equation*}

\section*{Appendix B (convergence rate)}
Let us show that the method \eqref{UniEBM} is first order accurate.
Toward that end we consider a time interval $(t_n, t_{n+1})$.
Let $\mathbf{C}_{\text{i}}(t_{n+1})$ be the exact solution to the problem \eqref{puba44}
with the initial condition $\mathbf C_{\text{i}}|_{t=t_n} = {}^{n} \mathbf{C}_{\text{i}}$.
It is sufficient to show that there exists $C < \infty$ such that
\begin{equation*}\label{AppendixB1}
\|{}^{n+1} \mathbf{C}_{\text{i}} - \mathbf{C}_{\text{i}}(t_{n+1}) \|
\leq C \ (\Delta t)^2 \quad \text{as} \ \Delta t \rightarrow 0.
\end{equation*}
It is well known that the EBM is first order accurate. Moreover, there
exists $\quad C_1 < \infty$ such that
\begin{equation*}\label{AppendixB2}
\|{}^{n+1} \mathbf{C}^{\text{EBM}}_{\text{i}} - \mathbf{C}_{\text{i}}(t_{n+1}) \|
\leq C_1 \ (\Delta t)^2 \quad \text{as} \ \Delta t \rightarrow 0.
\end{equation*}
On the other hand, $(\det(\mathbf{X}))^{-1/3}$ is a smooth function of $\mathbf{X}$
in vicinity of the exact solution $\mathbf{C}_{\text{i}}(t_{n+1})$.
Therefore, for sufficiently small time steps,
there exists a constant $C_3 < \infty$ such that
\begin{multline*}\label{AppendixB3}
|(\det({}^{n+1} \mathbf{C}^{\text{EBM}}_{\text{i}})^{-1/3} - 1 | =
|(\det({}^{n+1} \mathbf{C}^{\text{EBM}}_{\text{i}})^{-1/3} - \det(\mathbf{C}_{\text{i}}(t_{n+1}))^{-1/3} | \leq \\
  C_3 \ \|{}^{n+1} \mathbf{C}^{\text{EBM}}_{\text{i}} - \mathbf{C}_{\text{i}}(t_{n+1}) \|
\leq C_1 C_3 \ (\Delta t)^2 \quad \text{as} \ \Delta t \rightarrow 0.
\end{multline*}
Thus, for small $\Delta t$, there exists $C_4 < \infty$ such that
\begin{equation*}\label{AppendixB4}
\|\overline{{}^{n+1} \mathbf{C}^{\text{EBM}}_{\text{i}}} - {}^{n+1} \mathbf{C}^{\text{EBM}}_{\text{i}} \| =
|(\det({}^{n+1} \mathbf{C}^{\text{EBM}}_{\text{i}})^{-1/3} - 1 | \
\| {}^{n+1} \mathbf{C}^{\text{EBM}}_{\text{i}} \| \leq  C_4 \ (\Delta t)^2.
\end{equation*}
Finally, using the triangle inequality, we get the required estimation
\begin{multline*}\label{AppendixB5}
\|{}^{n+1} \mathbf{C}_{\text{i}} - \mathbf{C}_{\text{i}}(t_{n+1}) \| =
\|\overline{{}^{n+1} \mathbf{C}^{\text{EBM}}_{\text{i}}} - \mathbf{C}_{\text{i}}(t_{n+1}) \| \leq \\
\|\overline{{}^{n+1} \mathbf{C}^{\text{EBM}}_{\text{i}}} - {}^{n+1} \mathbf{C}^{\text{EBM}}_{\text{i}} \| +
\|{}^{n+1} \mathbf{C}^{\text{EBM}}_{\text{i}} - \mathbf{C}_{\text{i}}(t_{n+1}) \|
 \leq  (C_4 + C_1) \ (\Delta t)^2.
\end{multline*}

\section*{Appendix C (energy release during relaxation)}
Let us show that $\psi(\mathbf C \ {{}^{n+1}\mathbf C_{\text{i}}}^{-1})$
is a monotonically decreasing function of $\Delta t$ for a fixed $\mathbf C$,
where ${}^{n+1}\mathbf C_{\text{i}}$ is given by the explicit update formula \eqref{ExplUpdate}.
We make use of the fact that the free energy $\psi(\mathbf C \ {\mathbf C_{\text{i}}}^{-1})$ is
invariant under the isochoric change of the reference configuration.
Thus, due to the invariance of the
integration algorithm with respect to the reference change,
we may assume that $\overline{\mathbf C} = \mathbf{1}$.
In that case, the update formula takes the simple form
\begin{equation}\label{AppendixC1}
{}^{n+1} \mathbf{C}_{\text{i}}  =
\overline{\mathbf{\Phi}}, \quad
\mathbf{\Phi} := {}^{n} \mathbf{C}_{\text{i}} + \frac{\mu \Delta t}{\eta} \ \mathbf{1}.
\end{equation}
For $\overline{\mathbf C} = \mathbf{1}$, it follows from \eqref{freeen2} that
\begin{equation}\label{AppendixC2}
\frac{d}{d \Delta t} \psi(\mathbf C \ {}^{n+1}{\mathbf C_{\text{i}}}^{-1}) =
\frac{\mu}{2 \rho_{\scriptscriptstyle \text{R}} } \
\frac{d}{d \Delta t} \big( \text{tr} ( {\mathbf C_{\text{i}}}^{-1}) \big) =
- \frac{\mu}{2 \rho_{\scriptscriptstyle \text{R}} } \
\text{tr}\big( {}^{n+1}{\mathbf C_{\text{i}}}^{-1} \
\frac{d}{d \Delta t} ({}^{n+1} \mathbf{C}_{\text{i}})  \
{}^{n+1}{\mathbf C_{\text{i}}}^{-1} \big).
\end{equation}
Note that
\begin{equation}\label{AppendixC3}
\frac{d}{d \Delta t} \mathbf{\Phi} =  \frac{\mu}{\eta}
\mathbf{1}, \quad
\frac{d}{d \Delta t} (\det \mathbf{\Phi})^{-1/3} = -\frac{1}{3} \frac{\mu}{\eta} (\det \mathbf{\Phi})^{-1/3}
\ \text{tr}(\mathbf{\Phi}^{-1}).
\end{equation}
Thus, differentiating $(\ref{AppendixC1})_1$ and using \eqref{AppendixC3}, we get
\begin{equation*}\label{AppendixC4}
\frac{d}{d \Delta t} \ {}^{n+1} \mathbf{C}_{\text{i}}= (\det \mathbf{\Phi})^{-1/3}
\frac{\mu}{\eta} \big[\mathbf{1} - \frac{1}{3}\text{tr}(\mathbf{\Phi}^{-1}) \mathbf{\Phi}\big].
\end{equation*}
Substituting this result into  \eqref{AppendixC2} and taking into account
that ${}^{n+1}\mathbf C_{\text{i}}^{-1} = (\det \mathbf{\Phi})^{1/3} \mathbf{\Phi}^{-1}$, we get
\begin{multline*}\label{AppendixC2dss}
\frac{d}{d \Delta t} \psi(\mathbf C \ {}^{n+1}{\mathbf C_{\text{i}}}^{-1}) = -
(\det \mathbf{\Phi})^{1/3} \ \frac{\mu^2}{2 \rho_{\scriptscriptstyle \text{R}} \eta }  \
\text{tr}(\mathbf{\Phi}^{-1} [\mathbf{1} - \frac{1}{3}\text{tr}(\mathbf{\Phi}^{-1}) \mathbf{\Phi}] \mathbf{\Phi}^{-1} ) = \\
- (\det \mathbf{\Phi})^{1/3} \ \frac{\mu^2 }{2 \rho_{\scriptscriptstyle \text{R}} \eta } \text{tr}\big(\mathbf{\Phi}^{-1}
(\mathbf{\Phi}^{-1})^{\text{D}}\big) = - (\det \mathbf{\Phi})^{1/3} \ \frac{\mu^2}{2 \rho_{\scriptscriptstyle \text{R}} \eta }
\|(\mathbf{\Phi}^{-1})^{\text{D}} \|^2 \leq 0.
\end{multline*}
In particular, the free energy of the numerical solution for relaxation processes is a decreasing function of
$\Delta t$.

\section*{Appendix D (consistent tangent operator)}
Combining \eqref{2PKd} and \eqref{ExplUpdate}, we define a function $\tilde{\mathbf T}_1 (\mathbf C)$ through
\begin{equation}\label{AppendixD1}
\tilde{\mathbf T}_1 (\mathbf C) =  \mu \ \mathbf C^{-1}
\Bigg(\overline{\mathbf C} \Big(\overline{{}^{n} \mathbf{C}_{\text{i}} +  \frac{\mu \Delta t}{\eta}
\overline{\mathbf C}} \Big)^{-1}\Bigg)^{\text{D}}.
\end{equation}
Next, we abbreviate
\begin{equation}\label{AppendixD2}
\mathbf{\Phi} := {}^{n} \mathbf{C}_{\text{i}} +  \frac{\mu \Delta t}{\eta}
\overline{\mathbf C}, \quad
\mathbf{C}_{\text{i}}(\mathbf{C}) := \overline{{}^{n} \mathbf{C}_{\text{i}} +  \frac{\mu \Delta t}{\eta}
\overline{\mathbf C}} = \overline{\mathbf{\Phi}}.
\end{equation}
Using notations from \cite{ShutovKrKoo},
the consistent tangent operator is explicitly given by
\begin{multline}\label{AppendixD3}
\frac{\displaystyle \partial
\tilde{\mathbf T}_1 (\mathbf C)}{\displaystyle \partial \mathbf C} \stackrel{(\ref{AppendixD1})}{=}
\frac{\displaystyle \partial \big( \mu
\mathbf C^{-1} \ (\overline{\mathbf C}  \mathbf{C}_{\text{i}}^{-1})^{\text{D}} \big)}{\displaystyle \partial \mathbf C} =
- \mu \mathbf C^{-1} \odot \big(\mathbf C^{-1} (\overline{\mathbf C}  \mathbf{C}_{\text{i}}^{-1})^{\text{D}}\big)+ \\
\mu (\det\mathbf C)^{-1/3} \ \mathbf C^{-1} \cdot
\Big\{ - \frac{1}{3}(\mathbf C  \mathbf{C}_{\text{i}}^{-1})^{\text{D}} \ \otimes \mathbf C^{-1} +
\mathbb{P} : \big(\mathbb{I} \cdot \mathbf{C}_{\text{i}}^{-1} + \mathbf{C} \cdot
\frac{\displaystyle \partial
\mathbf{C}_{\text{i}}^{-1}}{\displaystyle \partial \mathbf C} \big) \Big\},
\end{multline}
\begin{equation}\label{AppendixD4}
\frac{\displaystyle \partial
\mathbf{C}_{\text{i}}^{-1}}{\displaystyle \partial \mathbf C} = -
(\mathbf{C}_{\text{i}}^{-1} \odot \mathbf{C}_{\text{i}}^{-1}): \frac{\displaystyle \partial
\mathbf{C}_{\text{i}}}{\displaystyle \partial \mathbf C},
\end{equation}
\begin{equation}\label{AppendixD5}
\frac{\displaystyle \partial
\mathbf{C}_{\text{i}}}{\displaystyle \partial \mathbf C} \stackrel{(\ref{AppendixD2})_2}{=}
(\det \mathbf \Phi)^{-1/3} \
\Big(\mathbb{I} - \frac{1}{3} (\mathbf \Phi \otimes \mathbf \Phi^{-1}) \Big) : \frac{\displaystyle \partial
\mathbf{\Phi}}{\displaystyle \partial \mathbf C},
\end{equation}
\begin{equation}\label{AppendixD6}
\frac{\displaystyle \partial
\mathbf{\Phi}}{\displaystyle \partial \mathbf C}  \stackrel{(\ref{AppendixD2})_1}{=}  \frac{\mu \Delta t}{\eta} (\det \mathbf C)^{-1/3} \
\Big(\mathbb{I} - \frac{1}{3} (\mathbf C \otimes \mathbf C^{-1}) \Big).
\end{equation}

Now let us prove the symmetry of the tangent operator $\frac{\displaystyle \partial
\tilde{\mathbf T}_1 (\mathbf C)}{\displaystyle \partial \mathbf C}$.
First, note that the material tangent resulting for frozen
inelastic flow ($\mathbf{C}_{\text{i}} = const$) corresponds to
the symmetric tangent operator of a hyperelastic material:
\begin{equation*}\label{AppendixD7}
\frac{\displaystyle \partial \big( \mu
\mathbf C^{-1} \ (\overline{\mathbf C}
\mathbf{C}_{\text{i}}^{-1})^{\text{D}} \big)}{\displaystyle
\partial \mathbf C}|_{\mathbf{C}_{\text{i}} = const}  \in Sym.
\end{equation*}
Thus, the symmetry of $\frac{\displaystyle \partial
\tilde{\mathbf T}_1 (\mathbf C)}{\displaystyle \partial \mathbf C}$ is equivalent
to the symmetry of $\mathbb{T}$ defined through
\begin{multline*}\label{AppendixD8}
\mathbb{T} := \frac{\displaystyle \partial
\tilde{\mathbf T}_1 (\mathbf C)}{\displaystyle \partial \mathbf C} -
\frac{\displaystyle \partial \big( \mu
\mathbf C^{-1} \ (\overline{\mathbf C}
\mathbf{C}_{\text{i}}^{-1})^{\text{D}} \big)}{\displaystyle \partial
\mathbf C}|_{\mathbf{C}_{\text{i}} = const} = \\
\mu (\det\mathbf C)^{-1/3} \ \mathbf C^{-1} \cdot
\Big\{\mathbb{P} : \big( \mathbf{C} \cdot
\frac{\displaystyle \partial
\mathbf{C}_{\text{i}}^{-1}}{\displaystyle \partial \mathbf C} \big) \Big\}.
\end{multline*}
It follows from \eqref{AppendixD4}--\eqref{AppendixD6} that for any $\mathbf{Y} \in Sym$
\begin{equation}\label{AppendixD9}
\frac{\displaystyle \partial
\mathbf{C}_{\text{i}}^{-1}}{\displaystyle \partial \mathbf C} : \mathbf{Y} = -  \frac{\mu \Delta t}{\eta}
\big(\mathbf{\Phi}^{-1}
(\mathbf{Y} \mathbf{C}^{-1})^{\text{D}} \overline{\mathbf{C}} \big)^{\text{D}} \mathbf{C}_{\text{i}}^{-1}.
\end{equation}
Note that
$\text{tr}(\mathbf{A} \mathbf{B}) =
\text{tr}(\mathbf{B} \mathbf{A})$ and
$\text{tr}(\mathbf{A} \mathbf{B}^{\text{D}}) = \text{tr}(\mathbf{A}^{\text{D}} \mathbf{B})
= \text{tr}(\mathbf{A}^{\text{D}} \mathbf{B}^{\text{D}})$ for arbitrary $\mathbf{A}, \mathbf{B}$.
Using these properties, we get for any $\mathbf{X}, \mathbf{Y} \in Sym$
\begin{multline*}\label{AppendixD10}
\mathbf{X} : \mathbb{T} : \mathbf{Y} =
\mu (\det\mathbf C)^{-1/3} \ \mathbf{X} : \Bigg( \mathbf C^{-1}  \Big( \mathbf C \
\frac{\displaystyle \partial
\mathbf{C}_{\text{i}}^{-1}}{\displaystyle \partial
\mathbf C} : \mathbf{Y} \Big)^{\text{D}} \Bigg) = \\
\mu (\det\mathbf C)^{-1/3} \text{tr}\Bigg\{ (\mathbf{X} \mathbf C^{-1})^{\text{D}}
\Big( \mathbf C \
\frac{\displaystyle \partial
\mathbf{C}_{\text{i}}^{-1}}{\displaystyle \partial \mathbf C} : \mathbf{Y} \Big)\Bigg\}
\stackrel{\eqref{AppendixD9}}{=}  \\
- \frac{\mu^2 \ \Delta t}{\eta}  \ (\det\mathbf C)^{-1/3} \
\text{tr}\Bigg\{ (\mathbf{X} \mathbf C^{-1})^{\text{D}}
\Big( \mathbf C \big(\mathbf{\Phi}^{-1}
(\mathbf{Y} \mathbf{C}^{-1})^{\text{D}} \overline{\mathbf{C}}
\big)^{\text{D}} \mathbf{C}_{\text{i}}^{-1} \Big)\Bigg\} = \\
- \frac{\mu^2 \ \Delta t}{\eta}   \ (\det\mathbf C)^{-1/3} \
\text{tr}\Bigg\{ \mathbf{C}_{\text{i}}^{-1}
(\mathbf{X} \mathbf C^{-1})^{\text{D}}
\mathbf C \big(\mathbf{\Phi}^{-1} \
(\mathbf{Y} \mathbf{C}^{-1})^{\text{D}} \overline{\mathbf{C}} \big)^{\text{D}} \Bigg\}  = \\
- \frac{\mu^2 \ \Delta t}{\eta}  \ (\det\mathbf C)^{-2/3} \  (\det\mathbf \Phi)^{1/3} \
\text{tr}\Bigg\{ \big( \mathbf{\Phi}^{-1} (\mathbf{X} \mathbf C^{-1})^{\text{D}}
\mathbf C \big)^{\text{D}} \ \big(\mathbf{\Phi}^{-1} \
(\mathbf{Y} \mathbf{C}^{-1})^{\text{D}} \mathbf{C} \big)^{\text{D}} \Bigg\}.
\end{multline*}
Obviously, this expression is symmetric
with respect to $\mathbf{X}$ and $\mathbf{Y}$.
Thus, $\mathbf{X} : \mathbb{T} : \mathbf{Y} = \mathbf{Y} : \mathbb{T} : \mathbf{X}$. The symmetry of the
consistent tangent $\frac{\displaystyle \partial
\tilde{\mathbf T}_1 (\mathbf C)}{\displaystyle \partial \mathbf C}$ is thus proved.

\section*{Appendix E (trajectory of the exact solution)}
Consider the initial value problem \eqref{puba44}, \eqref{initCond} with $\mathbf C(t) = const$.
Suppose that the numerical solution at time instance
$t_n$ is given by
\begin{equation}\label{AppendixE1}
{}^{n} \mathbf C_{\text{i}} = \overline{\mathbf C_{\text{i}}^0 +  {}^n \varphi
\overline{\mathbf C}}
\end{equation}
with some suitable ${}^n \varphi \geq 0$ . Substituting this into
the update formula \eqref{ExplUpdate} we get
\begin{equation*}\label{AppendixE2}
{}^{n+1} \mathbf{C}_{\text{i}}  =
\overline{{}^{n} \mathbf{C}_{\text{i}} +  \frac{\mu \Delta t}{\eta} \
\overline{\mathbf C}} =
\overline{\overline{\mathbf C_{\text{i}}^0 +  {}^n \varphi
\overline{\mathbf C}} +   \frac{\mu \Delta t}{\eta}
\overline{\mathbf C}} = \overline{\mathbf C_{\text{i}}^0 +  {}^{n+1} \varphi
\overline{\mathbf C}}
\end{equation*}
with a suitable ${}^{n+1} \varphi \geq 0$. Since the assumption \eqref{AppendixE1} is satisfied
for $n=1$, it is satisfied for any $n=1,2,3,...$ .
Finally, it is known that the numerical solution converges to the exact solution as $\Delta t \rightarrow 0$.
Therefore, for any $t \geq t^0$, we get from \eqref{AppendixE1}
\begin{equation}\label{AppendixE32}
\mathbf C_{\text{i}}(t) = \overline{\mathbf C_{\text{i}}^0 + \varphi(t)
\overline{\mathbf C}}.
\end{equation}
Thus, the trajectory of the exact solution is given by \eqref{AppendixE32}. By substituting
this relation into \eqref{puba44}, the following initial value problem
for $\varphi$ is obtained
\begin{equation*}\label{AppendixE3}
\dot{\varphi} =  \frac{\mu}{\eta} \big(\det(\mathbf C_{\text{i}}^0 + \varphi \overline{\mathbf C} )\big)^{1/3}, \quad
\varphi|_{t=t^0} =0.
\end{equation*}
In other words, the parametrization \eqref{AppendixE32} allows to reduce the six-dimensional flow rule to a one-dimensional one.

Finally, we have $\varphi \rightarrow 0$ as $t \rightarrow t^0$. Thus, since
$\det(\mathbf C_{\text{i}}^0)=1$ we get  $\det(\mathbf C_{\text{i}}^0 + \varphi(t) \overline{\mathbf C} ) \approx 1$
as $t \rightarrow t^0$. Substituting this into \eqref{AppendixE3}, we get $\varphi(t) \approx  \frac{\mu (t-t^0)}{\eta}$
as $t \rightarrow t^0$.


\begin{thebibliography}{10}

\bibitem{BalanTsakmakis}{C. Balan, C. Tsakmakis,
A finite deformation formulation of the 3-parameter viscoelastic fluid,
Journal of non-newtonian fluid mechanics, 103(1) (2002), 45--64.}

\bibitem{Broec}
C. Br\"ocker, A. Matzenmiller,
An enhanced concept of rheological models to represent
nonlinear thermoviscoplasticity and its energy storage
behavior, Continuum Mech. Thermodyn.
DOI 10.1007/s00161-012-0268-3

\bibitem{DettRes}{W. Dettmer, S. Reese,  On the theoretical and numerical modelling of
Armstrong–-Frederick kinematic hardening in the finite strain regime,
Computer Methods in Applied Mechanics and Engineering, 193 (2004) 87–-116.}

\bibitem{Drozdov}{A. Drozdov, Viscoelastic Structures:
Mechanics of Growth and Aging, Academic Press, 1998}

\bibitem{Eidel}{B. Eidel, C. Kuhn, Order reduction in computational inelasticity: Why it happens and
how to overcome it -- The ODE-case of viscoelasticity, International
Journal for Numerical Methods in Engineering, 87 (2011) 1046--1073.}

\bibitem{EidelStumpf}{B. Eidel, F. Stumpf, J. Schr\"oder, 
Finite strain viscoelasticity: how to consistently couple discretizations in time and space
on quadrature-point level for full order $p \geq 2$ and a considerable speed-up,
Computational Mechanics, (2013) DOI 10.1007/s00466-013-0823-1}

\bibitem{Feigen}{H. P. Feigenbaum, J. Dugdale, Y. F. Dafalias, K. I. Kourousis, J. Plesek,
Multiaxial ratcheting with advanced kinematic and directional distortional
hardening rules, International Journal of Solids and Structures, 49 (2012) 3063--3076.}

\bibitem{GassFor}{T. C. Gasser, C. Forsell, The numerical implementation of invariant-based viscoelastic
formulations at finite strains. An anisotropic model for the passive myocardium,
Computer Methods in Applied Mechanics and Engineering, 200 (2011) 3637--3645.}

\bibitem{Gent 1960}{A. N. Gent, Simple rotary dynamic testing machine.
Rubber Chemistry and Technology, 34(3) (1961) 790--794.}

\bibitem{GovRes}{S. Govindjee, S. Reese, A Presentation and
comparison of two large deformation viscoelasticity models,
Journal of Engineering Materials and Technology, 119, (1997) 251--255.}

\bibitem{Hartmann}{S. Hartmann, G. L\"uhrs, P. Haupt, An efficient
stress algorithm with applications in viscoplasticity and plasticity, International Journal
for Numerical Methods in Engineering, 40 (1997) 991--1013.}

\bibitem{Hartmann2002}{S. Hartmann, Computation in finite-strain viscoelasticity:
finite elements based on the interpretation as differential-algebraic equations,
Computer Methods in Applied Mechanics and Engineering, 191 (2002) 1439--1470.}

\bibitem{HartmannHabil}{S. Hartmann, Finite-Elemente Berechnung inelastischer Kontinua.
Interpretation als Algebro-Differentialgleichungssysteme, Habilitation thesis, Kassel, 2003.}

\bibitem{Hartmann2010}{S. Hartmann, A.-W. Hamkar, Rosenbrock-type methods
applied to finite element computations within finite
strain viscoelasticity,
Computer Methods in Applied Mechanics and Engineering, 199 (2010) 1455--1470.}

\bibitem{Hasanpour}{K. Hasanpour, S. Ziaei-Rad, M. Mahzoon, A large deformation
framework for compressible viscoelastic materials: Constitutive equations and finite
element implementation, International Journal of Plasticity, 25 (2009) 1154--1176.}

\bibitem{Haupt}{P. Haupt, Continuum Mechanics and Theory of Materials, 2nd edition,
Springer, 2002.}

\bibitem{HauptLi}{P. Haupt, A. Lion, On finite linear viscoelasticity
of incompressible isotropic materials, Acta Mechanica 159, 87-124, 2002.}

\bibitem{Helm1}{D. Helm,   Formged\"achtnislegierungen, experimentelle Untersuchung,
ph\"anomenologische Modellierung und numerische Simulation der thermomechanischen
Materialeigenschaften, Universit\"atsbibliothek Kassel, 2001.}

\bibitem{Helm2}{D. Helm, Stress computation in finite
thermoviscoplasticity. International Journal of Plasticity, 22 (2006) 1699--1721.}

\bibitem{Helm3}{D. Helm, Thermomechanics of martensitic phase transitions in shape
memory alloys I, constitutive theories for small and large deformations, J. Mech.
Mater. Struct., 2(1) (2007)  87--112.}

\bibitem{Hohl 2007}{C. Hohl. Anwendung der Finite-Elemente-Methode zur Parameteridentifikation
und Bauteilsimulation bei Elastomeren mit Mullins-Effekt. D\"usseldorf: VDI Verlag 2007}

\bibitem{HolLoug}{D. W. Holmes, J. G. Loughran, Numerical
aspects associated with the implementation of a finite
strain, elasto-viscoelastic-viscoplastic constitutive theory in
principal stretches, Int. J. Numer. Meth. Engng, 83 (2010) 366--402.}

\bibitem{Holzapfel}
{G.A. Holzapfel, On large strain viscoelasticity: continuum formulation and finite element applications to
elastomeric structures. Int. J. Numer. Meth. Eng. 39, (1996) 3903--3926.}

\bibitem{HubTsak}{N. Huber, C. Tsakmakis, Finite deformation viscoelasticity laws, Mechanics of Materials 32 (2000) 1--18}

\bibitem{Itskov}
M.~Itskov, Tensor Algebra and Tensor Analysis for Engineers:
With Applications to Continuum Mechanics
(Springer, 2007).

\bibitem{Johlitz}{M. Johlitz, H. Steeb, S. Diebels,
A. Chatzouridou, J. Batal, W. Possart, Experimental and theoretical investigation of nonlinear
viscoelastic polyurethane systems, J. Mater. Sci., 42 (2007) 9894--9904}

\bibitem{JohnsonSegalman}{M. W. Johnson, D. Segalman, A model for viscoelastic fluid behavior
which allows non-affine deformation. J. Non-Newtonian Fluid Mech. 2,
(1977) 255--270.}

\bibitem{Klauke1}{R. Klauke, T. Alshuth, J. Ihlemann,
 Lifetime prediction of rubber products under simple-shear loads with rotary axes.
 In: G. Heinrich, M. Kaliske, A. Lion, S. Reese (Editors):
 Constitutive Models for Rubber VI, Taylor \& Francis Group, London, (2009) 235--240.}

\bibitem{Klauke2}{R. Klauke, T. Alshuth, J. Ihlemann,
Lebensdauervorhersage von technischen Gummiwerkstoffen
unter einfacher Scherung mit rotierenden Achsen.
Kautschuk Gummi Kunststoffe 63, (2010) 286--290.}

\bibitem{Kleuter}{B. Kleuter, A. Menzel, P. Steinmann,
Generalized parameter identification for finite viscoelasticity,
Computer Methods in Applied Mechanics and Engineering, 196 (2007) 3315--3334.}

\bibitem{Koprowski}{N. Koprowski-Thei\ss, M. Johlitz, S. Diebels,
Modelling of a cellular rubber with nonlinear viscosity
functions, Experimental Mechanics, 51 (2011) 749--765.}

\bibitem{Kroener}
E. Kr\"oner, Allgemeine Kontinuumstheorie der Versetzungen und
Eigenspannungen, Arch. Rational Mech. Anal., 4 (1959) 273--334.

\bibitem{LaMantia}{F. P. La Mantia, Non linear viscoelasticity of polymeric liquids interpreted by
means of a stress dependence of free volume. Rheol. Acta. 16, 302--308 (1977)}

\bibitem{Landgraf} R. Landgraf, J. Ihlemann, Vergleich zweier Ans\"atze zur Beschreibung
nichtlinearer Viskoelastizit\"at auf Basis des Maxwell-Elements, PAMM, (10) 1, 303--304.

\bibitem{Lee}
E. H. Lee, Elastic–-plastic deformation at finite strains,
J. Appl. Mech., 36 (1969) 1--6.

\bibitem{Lejeunes}{S. Lejeunes, A. Boukamel, S. M\'eo,
Finite element implementation of nearly-incompressible rheological models
based on multiplicative decompositions, Computers and Structures 89, (2011) 411--421.}

\bibitem{Leonov}
A. I. Leonov, Nonequilibrium thermodynamics and rheology of viscoelastic
polymeric media. Rheol. Acta 15,  (1976) 85--98.

\bibitem{LionAM}{A. Lion, A physically based method to represent
the thermo-mechanical behaviour of elastomers, Acta Mechanica 123, 1-25 (1997)}

\bibitem{LionIJP}{A. Lion, Constitutive modelling in finite thermoviscoplasticity:
a physical approach based on nonlinear rheological elements,
International Journal of Plasticity, 16 (2000) 469--494. }

\bibitem{LionHab}{A. Lion, Thermomechanik von Elastomeren, Habilitation  thesis, Kassel, 2000.}

\bibitem{Meng}{X. N. Meng, T. A. Laursen, Energy consistent algorithms for
dynamic finite deformation plasticity,
Computer Methods in Applied Mechanics and Engineering, 191 (2002) 1639--1675.}

\bibitem{MARC}{MSC.Software Corporation: MSC.Marc 2010,  Volume A: Theory and User Information.}

\bibitem{Nedjar}{B. Nedjar, Frameworks for finite strain viscoelastic-plasticity based on
multiplicative decompositions. Part I: Continuum formulations, Comput. Methods Appl. Mech.
Engrg. 191 (2002) 1541--1562.}

\bibitem{Nedjar2}{B. Nedjar, Frameworks for finite strain viscoelastic-plasticity based on
multiplicative decompositions. Part II:
Computational aspects, Comput. Methods Appl. Mech.
Engrg. 191 (2002) 1563--1593.}

\bibitem{Neff}
P. Neff, Mathematische Analyse multiplikativer Viskoplastizit{\"a}t.
Ph.D. Thesis TU Darmstadt.
(Shaker Verlag, 2000).

\bibitem{NishiguchiA}
I. Nishiguchi, T.-L. Sham, E. Krempl, A finite deformation theory
of viscoplasticity based on overstress: Part I-Constitutive Equations. Trans.
ASME J Appl. Mech. 57, (1990) 548--552.

\bibitem{NishiguchiB} I. Nishiguchi, T.-L. Sham, E. Krempl, A finite deformation theory
of viscoplasticity based on overstress: Part ll-Finite element implementation and
numerical experiments. Trans. ASME J Appl. Mech. 57, (1990) 553--561.

\bibitem{PearsonMiddle}{G. Pearson, S. Middleman, Elongation flow behavior of viscoelastic
liquids: modelling bubble dynamics with viscoelastic constitutive relations.
Rheol. Acta 17, (1978) 500--510.}

\bibitem{PericCr}{D. Peric, A.J.L. Crook, Computational strategies for predictive geology with
reference to salt tectonics, Comput. Methods Appl. Mech. Engrg. 193 (2004) 5195--5222}

\bibitem{Rauchs}{G. Rauchs, Finite element implementation including sensitivity analysis of a simple finite
strain viscoelastic constitutive law, Computers and Structures 88 (2010) 825--836}

\bibitem{ReeseG}{S. Reese, S. Govindjee, A theory of finite viscoelasticity and numerical aspects,
International Journal of Solids and Structures, 35 (1998), 3455--3482.}

\bibitem{ReeseHab}{S. Reese, Thermomechanische Modellierung gummiartiger Polymerstrukturen,
Habilitation thesis, Hannover, 2000.}

\bibitem{Reiner}{M. Reiner, Deformation, Strain and Flow. An Elementary Introduction to Rheology, 2nd edition, 1960.}

\bibitem{ShutovKrVisc}
A. V. Shutov, R. Krei{\ss}ig, Finite strain viscoplasticity with
nonlinear kinematic hardening: Phenomenological modeling and
time integration, Computer Methods in Applied
Mechanics and Engineering, \textbf{197}, 2015--2029 (2008).

\bibitem{ShutovKrKoo}{A. V. Shutov, R.~Krei{\ss}ig, Application
of a coordinate-free tensor formalism to the numerical
implementation of a material model, ZAMM, \textbf{88}, 11, 888-909 (2008).}

\bibitem{ShutovKrStab}{A. V. Shutov, R. Krei{\ss}ig,
Geometric integrators for multiplicative viscoplasticity: Analysis
of error accumulation, Comput. Methods Appl. Mech. Engrg. 199 (2010) 700--711.}

\bibitem{ShutovKuprin}{A. V. Shutov, C. Kuprin, J. Ihlemann, M. F.-X.Wagner, C. Silbermann,
Experimentelle Untersuchung und numerische Simulation
des inkrementellen Umformverhaltens von Stahl 42CrMo4, Mat.-wiss. u.Werkstofftech., 41(9), (2010) 765--775.}

\bibitem{ShutovPaKr}{A. V. Shutov, S. Panhans, R. Krei{\ss}ig, A phenomenological model
of finite strain viscoplasticity
with distortional hardening, ZAMM, \textbf{91}, 8, 653-680 (2011).}

\bibitem{ShutovPfeIhl}{A. V. Shutov, S. Pfeiffer, J. Ihlemann,
On the simulation of multi-stage forming processes:
invariance under change of the reference configuration,
Mat.-wiss. u.Werkstofftech., 43(7), (2012) 617--625.}

\bibitem{ShutovIhle}{A. V. Shutov, J.~Ihlemann, A viscoplasticity model with an enhanced control of the yield
surface distortion, International Journal of Plasticity, \textbf{39}, 152-167 (2012).}


\bibitem{SimoMeschke}{J. C. Simo, G. Meschke,
A new class of algorithms for classical plasticity
extended to finite strains. Application to geomaterials,
Computational mechanics, 11(4), (1993) 253--278.}


\bibitem{SimMieh} J. C. Simo, C. Miehe, Associative coupled thermoplasticity at finite strains: formulation, numerical
analysis and implementation. Computer Methods in Applied Mechanics and Engineering 98, (1992) 41--104.

\bibitem{Simo} J. C. Simo, Algorithms for static and dynamic multiplicative plasticity that preserve the classical return
mapping schemes of the infinitesimal theory. Computer Methods in Applied Mechanics and Engineering 99, (1992) 61--112.

\bibitem{SimHug}{J. Simo, T. Hughes, Computational inelasticity, Springer, 1998.}

\bibitem{SussmanBathe1987} T. Sussman, K. J. Bathe, A finite element formulation for nonlinear
incompressible elastic and inelastic analysis. Computers and Structures, 26(1), (1987) 357--409.

\bibitem{Vladimirov}{I. Vladimirov, M. Pietryga, S. Reese, On the modelling of non-linear
kinematic hardening at finite strains with application to springback -- Comparison of time
integration algorithms, Int. J. Numer. Meth. Engng 75 (2008), 1--28.}

\bibitem{Wiechert}{E. Wiechert, "Ueber elastische Nachwirkung",
Dissertation, K\"onigsberg University, Germany, 1889.}

\bibitem{Yeoh 1993}{O. H. Yeoh, Some forms of the strain energy function for rubber,
Rubber Chemistry and technology, 66(5) (1993), 754-771.}

\end{thebibliography}
\end{document}